\documentclass[trsc,nonblindrev]{informs3} 

\OneAndAHalfSpacedXI 

\usepackage{natbib}
 \bibpunct[, ]{(}{)}{,}{a}{}{,}%
 \def\bibfont{\small}%
 \def\bibsep{\smallskipamount}%
\renewcommand{\bibfont}{\footnotesize}

 \usepackage[ruled, vlined, linesnumbered]{algorithm2e}
\TheoremsNumberedThrough     
\EquationsNumberedThrough    

\usepackage[ruled, vlined, linesnumbered]{algorithm2e}

\usepackage{algpseudocode}
\usepackage{amsmath}
\usepackage{amssymb}
\usepackage{tabularx}
\usepackage{booktabs}
\usepackage{multirow}
\usepackage{hyperref}
\usepackage{setspace}
\usepackage[ruled,vlined]{algorithm2e}
\usepackage{threeparttable}
\usepackage{subfigure}
\usepackage{bm}
\usepackage{float}
\usepackage{bbding}
\usepackage{tikz}
\usetikzlibrary{positioning}
\usetikzlibrary{shapes,arrows,arrows,positioning,fit}
\usepackage[utf8]{inputenc}
\usepackage[T1]{fontenc}
\usepackage{pgfplots}
\usepackage{pgfplotstable}
\usepackage{enumitem} 
\usepackage{relsize}
\usepackage[normalem]{ulem}
\usetikzlibrary{%
  arrows,%
  calc,
  shapes,
  arrows,
  shapes.misc,
  shapes.arrows,%
  chains,%
  matrix,%
  positioning,
  scopes,%
  decorations.pathmorphing,
  shadows%
}

\usepackage{xspace}
\newcommand{\oxc}{,\xspace} 
\newcommand{\SetCond}[2]{
  \{ #1 \mid #2 \}%
}

\begin{document}

\RUNAUTHOR{}

\RUNTITLE{}

\TITLE{\textbf{Benders Decomposition for Passenger-Oriented Train Timetabling with Hybrid Periodicity}}

\ARTICLEAUTHORS{%
\AUTHOR{Zhiyuan Yao}
\AFF{Frontiers Science Center for Smart High-Speed Railway System, Beijing Jiaotong University, 100044 Beijing, China \\ RPTU Kaiserslautern-Landau, 67663 Kaiserslautern, Germany\\
\EMAIL{zhiyuan@bjtu.edu.cn}}

\AUTHOR{Anita Schöbel}
\AFF{RPTU Kaiserslautern-Landau, 67663 Kaiserslautern, Germany \\ Fraunhofer Institute for Industrial Mathematics, 67663 Kaiserslautern, Germany\\
\EMAIL{anita.schoebel@rptu.de}}

\AUTHOR{Lei Nie}
\AFF{Frontiers Science Center for Smart High-Speed Railway System, Beijing Jiaotong University, 100044 Beijing, China \\ School of Traffic and Transportation, Lanzhou Jiaotong University, 730070 Lanzhou, China\\
\EMAIL{lnie@bjtu.edu.cn}}

\AUTHOR{Sven Jäger}
\AFF{RPTU Kaiserslautern-Landau, 67663 Kaiserslautern, Germany\\
TU Berlin, 10623 Berlin, Germany\\
\EMAIL{jaeger@math.tu-berlin.de}}
}

\ABSTRACT{
	Periodic timetables are widely adopted in passenger railway operations due to their regular service patterns and well-coordinated train connections. However, fluctuations in passenger demand require varying train services across different periods, necessitating adjustments to the periodic timetable. This study addresses a hybrid periodic train timetabling problem, which enhances the flexibility and demand responsiveness of a given periodic timetable through schedule adjustments and aperiodic train insertions, taking into account the rolling stock circulation. Since timetable modifications may affect initial passenger routes, passenger routing is incorporated into the problem to guide planning decisions towards a passenger-oriented objective. Using a time-space network representation, the problem is formulated as a dynamic railway service network design model with resource constraints. To handle the complexity of real-world instances, we propose a decomposition-based algorithm integrating Benders decomposition and column generation, enhanced with multiple preprocessing and accelerating techniques. Numerical experiments demonstrate the effectiveness of the algorithm and highlight the advantage of hybrid periodic timetables in reducing passenger travel costs.
}


\maketitle

\linespread{1.38}\selectfont

\section{Introduction}
	In passenger railway service planning, train timetables determine the arrival and departure times of trains, thus directly influencing the available train services that passengers can choose within a specific time frame. Periodic timetables, where trains operate on a repeated schedule with a fixed time interval (i.e., the cycle time), are widely implemented in European countries such as the Netherlands, Germany\oxc and Switzerland \citep{Peeters2003-thesis}. Due to their service regularity and well-coordinated train connections, periodic timetables are not only easy for passengers to remember, but can also attract more passengers and enhance the revenue of railway companies \citep{Johnson2006-RegularityBenefit}.
	
	The implementation of periodic timetables involves strategic and tactical planning processes carried out in several sequential stages \citep{Schiewe2022-Seq}. First, based on the hourly passenger demand, line planning is conducted to develop a line plan specifying the routes, frequencies and intermediate stops of the lines to be operated. This is followed by the train timetabling stage, where a periodic timetable for one cycle time is developed, specifying the arrival and departure times of individual trains at stations along their routes. Then, this periodic timetable is replicated and extended over the day to generate the trips, based on which the rolling stock circulation is determined. 
	
	The above planning stages produce an executable periodic timetable with regular train services throughout a day. However, since passenger demand is temporally different, such timetables may result in capacity wastes during off-peak hours and insufficient capacity during peak ones. It is also possible that certain train services or stop patterns are more desirable at specific times. Planning train operations for different periods simultaneously is hard, first approaches have been suggested in \cite{SchieweA2023-LPP-MultiPeriod}. The practical approach is to mitigate these issues using \emph{post-processing} strategies, such as providing heterogeneous rolling stock units (RSUs) with the suitable number of carriages for different periods \citep{Alfieri2006-UnitCompo}, and adopting hybrid periodic timetables where aperiodic train services are added to accommodate irregular passenger demand and several periodic services can be canceled \citep{Robenek2017-Hybrid}. From the modeling perspective, the first strategy concerns the efficient circulation of heterogeneous RSUs to match the demand, while the second one involves adjustments to the line plan and timetable across the day.
	
	Given the interdependence between these two strategies, this study seeks to integrate them to develop timetables with hybrid periodicity. From an operational perspective, the implementation of the line plan and timetable relies on the availability of RSUs --- scheduling an aperiodic train may require canceling a periodic train if no RSUs are available. Conversely, the line plan and timetable influence the connection opportunities between train services and thus impact rolling stock circulation. From the passenger's perspective, the line plan, timetable, and rolling stock circulation collectively determine passenger traveling options. For instance, if high-capacity RSUs are unavailable during peak hours, aperiodic trains with appropriate intermediate stops should be inserted to match the demand. Therefore, the line plan, timetable\oxc and rolling stock circulation shall be jointly optimized to offer a feasible and passenger-oriented railway service. 
	
	Existing studies have typically applied only one of the above two strategies, or not given much attention to passenger preferences. \cite{Tan2020-Extra} and \cite{Tan2020-Insert-Unit} inserted aperiodic train services into a periodic timetable, and the latter study included the rolling stock circulation of aperiodic train services. However, passenger routing was not included to guide the timetable adjustments, and the line plan of inserted trains was given beforehand. For studies integrating passenger perspectives, \cite{Li2023-PeriodicExpansion} integrated rolling stock circulation and passenger assignment for a \emph{fixed} one-day periodic timetable. \cite{Robenek2017-Hybrid} developed hybrid timetables with passenger routing, while rolling stock circulation and track capacity constraints were not taken into account. Distinct from the above studies, this paper aims to improve a given periodic timetable by simultaneously incorporating periodic timetable adjustments, line planning and timetabling of aperiodic trains, and rolling stock circulation, while considering practical requirements such as seat and track capacities. Passenger routing is further integrated to guide these decisions towards a passenger-oriented objective.
	
	The problem we consider is inherently large-scale, as the considered time horizon is multiple periods rather than a single one, which significantly increases the number of trains and time-dependent passenger groups. Additionally, the problem is highly complex, as it integrates line planning, train timetabling\oxc and RSU circulation, each of which is already challenging on its own. Moreover, since original passenger transfers in the periodic timetable might be influenced due to timetable adjustments, transfers must be included in passenger routing to take these impacts into account, which further increases the number of passenger paths. To solve such a complex problem, we develop a time-space network representation and formulate a dynamic railway service network design model with resource constraints. To tackle real-life instances, we propose an efficient decomposition-based algorithm combining Benders decomposition and column generation.
	
	Based on the above discussions, the contributions of this study are three-fold:
	\begin{itemize}
		\item A hybrid periodic train timetabling problem (HPTTP) is introduced, which improves a given periodic timetable by timetable adjustment, aperiodic train insertion\oxc and rolling stock circulation. Various planning decisions are guided by the dynamic distribution of passenger flows to increase the flexibility of a periodic train timetable in a passenger-oriented way.
		\item A dynamic railway service network design model with resource constraints is formulated, which is solved by a decomposition-based algorithm combining Benders decomposition and column generation. The Benders decomposition framework separates the overall problem into train service network design and passenger routing problems, while the column generation specifically addresses the large-scale passenger routing subproblem with transfer requirements. We also propose a problem variant that fixes routes for partial passenger groups, which performs well in large instances.
		\item Extensive numerical experiments are conducted on the artificial and real-world instances, demonstrating the effectiveness of the proposed solution approach, as well as the benefits of the hybrid periodicity in terms of reducing passenger travel costs.
	\end{itemize}

	The remainder of the paper is organized as follows. Section \ref{sec:lit rev} reviews existing studies related to HPTTP and highlights the novelty of our work. Section \ref{sec:problem into} describes the key characteristics of the studied problem. The time-space network representation for the problem is detailed in Section \ref{sec:ts net}, which serves as a basis for the arc-path dynamic service network design model formulated in Section \ref{sec:model}. To cope with challenging instances, an exact Benders decomposition algorithm is developed in Section \ref{sec:algo}. Section \ref{sec: subset routing} introduces a problem variant that fixes the routes of partial passenger groups, which achieves better results in large instances. Numerical experiments validating the effectiveness of the proposed methods are reported in Section \ref{sec:numerical}, and conclusions are drawn in Section \ref{sec:conclu}.
	
	\section{Literature review} \label{sec:lit rev}
	Although the HPTTP is a newly defined problem, it is closely related to the passenger-oriented line planning problem (LPP), train timetabling/rescheduling problem\oxc and rolling stock circulation problem. To keep the literature review brief and clear, we mainly review the integrated versions of these problems: the integrated line planning and timetabling problem, the train (re)scheduling problem with rolling stock and passenger considerations\oxc and the extra train scheduling problem. The following sections summarize the most relevant studies.
	
	\textbf{Integrated line planning and timetabling problem.}
	In HPTTP, the line plans and timetables of extra aperiodic train services are simultaneously optimized. Due to the strong interdependence between the two problems, their integration has received much attention in the literature. Iterative solution frameworks, where the line plan and timetable are alternately optimized, have been put forward \citep{Anita2017-Eigen, Burggraeve2017-LPP-TPP, YanF2019-LPP,Fuchs2022-LPP-TPP}. In these studies, when timetables are not feasible or good enough, certain modifications to the line plan are implemented, such as resolving conflicts, lowering maximum frequencies, increasing seat capacities or adjusting train stops. By encoding both line plan and timetable information, \cite{Kaspi2013-TTP-LPP-PR} directly addressed the integrated problem using a cross-entropy heuristic. Besides the heuristic approaches, integrated models haven also been proposed. \cite{ZhangCT2021-LPP} included line planning decisions in the event-activity network-based timetabling model, and \cite{ZhangYX2022-lntegrated} integrated a line planning model with a time-space network-based timetabling model. However, passenger transfers were not considered in these two studies.
	
	It is noted that passenger transfers have been considered in an easier problem setting, where only train stop decisions are considered in LPP. By expanding the time-space network to include nodes that imply train stopping decisions, the integrated stop planning and timetabling problem has been modeled as a dynamic service network design problem and solved by dual decomposition approaches \citep{Xu2021-LR, Yao2023-ADMM}. This modeling technique is straightforward; therefore we extend this approach to further include the origin-destination (OD) decisions of the extra trains.
	
	\textbf{Train (re)scheduling with passenger and rolling stock  considerations.}
	In HPTTP, the timetable is optimized in a passenger-oriented way, either by scheduling new aperiodic trains or adjusting existing periodic timetables. This makes the problem related to both the passenger-oriented train scheduling and rescheduling problems. \cite{Binder2017-TriObj} proposed an arc-based time-space network model for the multi-objective train rescheduling problem, integrating train cancellation, train rerouting, emergency train scheduling\oxc and passenger routing. Using a similar modeling technique, \cite{Zhan2021-ADMM} and \cite{Xu2021-LR} modeled passenger-oriented train rescheduling and timetabling problems that consider passenger transfers. They developed ADMM- and Lagrangian relaxation-based algorithms to solve large-scale instances, respectively. \cite{Zhu2020-ReschPax-EAN} proposed an event-activity network-based mixed-integer programming (MIP) model for timetable rescheduling with passenger reassignment to offer more flexible passenger routing alternatives during disruptions. \cite{Zhang2023-LR-Rer} proposed a Lagrangian heuristic for the train rescheduling problem under large-scale disruptions, aiming to minimize passenger delays.

	Further, rolling stock circulation has been integrated in the problem to specify the seat capacities of train services and ensure the availability of RSUs. Generally, there are two ways to take into account rolling stock resources within train scheduling. The first approach explicitly models complete rolling stock scheduling decisions, specifying the sequence of trips assigned to each RSU. \cite{Xu2018-TimeRoll} designed a time-space-state network for the integrated timetabling and locomotive assignment problem and proposed a locomotive flow-based formulation where locomotive paths explicitly represented trip sequences. Using a similar modeling technique, \cite{Wang2024-TTP-RSP} developed a branch-and-price algorithm for the integrated timetabling problem with a passenger-oriented objective. Besides the integrated approaches, \cite{Veelenturf2017-PaxDisrup} proposed an iterative framework for train rescheduling, passenger routing\oxc and rolling stock rescheduling, where simulated passenger behaviors guide the rescheduling decisions of timetables and RSUs. In the context of periodic timetabling, \cite{LiebMohr2007-PESP} incorporated vehicle scheduling into the periodic event scheduling problem (PESP) model and analyzed the number of vehicles needed to operate a periodic timetable under specific scenarios.
	
	The second approach only takes into account rolling stock balances at stations. \cite{Louwerse2014-MILP-Inventory} incorporated the rolling stock inventory events and activities into an event-activity network of the train rescheduling problem to maintain rolling stock balance at stations. This technique was then adopted by \cite{Liao2021-LR} in the time-space network representation of the problem. In a passenger-oriented context, \cite{Pan2023-TTP-RSP} further introduced flexible coupling and decoupling of RSUs in the time-space network representation. However, it is important to note that this approach does not determine the trip sequences of units, and it is possible that the resulting solution becomes infeasible if rolling stock maintenance restrictions are introduced. 
	
	\textbf{Extra train scheduling problem.}
	HPTTP is closely related to the extra train scheduling problem (ETSP), as aperiodic train services can be added to the basic periodic timetable. Considering the interactions between existing and extra train services, \cite{Burdett2009-ATSP} formulated ETSP as a hybrid job shop scheduling problem with general time window constraints. \cite{Cacchiani2010-LR-Extra} developed a Lagrangian heuristic to schedule extra freight trains in a railway network where passenger train timetables remain unchanged. \cite{Gao2018-ATSP} proposed a three-stage optimization method for scheduling extra trains in a rail corridor and minimized the extra trains' travel times and adjustments to the existing timetable. \cite{Tan2015-Doctoral} and \cite{Tan2020-Insert-Unit} formulated event-activity network-based models to insert extra aperiodic trains into a periodic timetable, and examined objectives such as minimizing adjustments to the initial periodic timetable and the number of RSUs required to operate the extra trains. \cite{Tan2020-Extra} furthered minimized the number of missed initial train connections in a similar problem setting. \cite{Dekker2024-InsertPath} proposed a dynamic programming approach to quickly insert new train paths into a fixed timetable. In general, ETSP extends the train timetabling problem by introducing strict/soft time windows for train services and places greater emphasis on the interaction between existing and extra trains instead of passengers' perspectives. Besides, line plans of extra trains are usually given.
	
	\textbf{Integrated problem and literature gaps.}
	In the hybrid periodic train timetabling problem, existing studies have typically focused on the feasibility of inserting aperiodic trains while paying little attention to the associated impacts on passengers \citep{Tan2015-Doctoral, Tan2020-Extra, Tan2020-Insert-Unit}, fixed the timetable without considering potential new train services \citep{Li2023-PeriodicExpansion}, or neglected critical track capacity constraints \citep{Robenek2017-Hybrid}. In contrast, our new version of the problem is able to systematically improve passenger service levels by integrating various operational decisions and passenger routing, as well as to enhance its applicability by considering track capacity and seat capacity constraints.
	
	Comparing with the train rescheduling literature, \cite{Binder2017-TriObj} proposed a similar problem formulation, while rolling stock circulation and line planning decisions weren't considered, and the proposed model could only solve relatively small instances with 2 hours of horizon and 24 trains. Larger instances were tackled by \cite{Veelenturf2017-PaxDisrup} with an iterative heuristic framework that alternately optimized relevant subproblems. However, the major difference between HPTTP and the rescheduling problem is that HPTTP is a tactical-level problem, which favors better quality of the results even though more computation times are needed. Therefore, an exact approach is needed for HPTTP to guarantee the quality of the results. In addition, unlike the limited disruption durations in rescheduling problems, HPTTP must address the planning horizon of multiple cycle periods and consider  larger numbers of trains and time-dependent passenger groups, greatly increasing its complexity.

	From the modeling perspective, HPTTP can be viewed as a train timetabling problem integrated with rolling stock circulation, line planning\oxc and passenger routing. \cite{Anita2017-Eigen} proposed an eigenmodel for iterating among these components. However, the framework was generic and did not include specific solution approaches. \cite{Schiewe2020-book} and \cite{Schiewe2022-Seq} formulated integrated models with coupling constraints across different systems. However, the model can only be solved in very small instances owing to the extreme complexity of the problem. 

	Based on the above discussions, we highlight the importance of the proposed efficient and exact Benders decomposition algorithm combined with column generation for such a complex problem. Various preprocessing and accelerating techniques are incorporated to speed up computation. The main framework of the algorithm not only suits HPTTP, but also applies to relevant rescheduling and timetabling problems.
	\section{Problem description} \label{sec:problem into}
	
	This section details the characteristics of the relevant components of HPTTP, and the notations used are summarized in Appendix~\ref{app: tables}. 
	
	\textbf{Line planning and timetabling.}
	The \emph{railway network} considered contains a set of stations $\cal M$ and a set of direction-dependent railway track sections ${\cal E} \subseteq {\cal M} \times {\cal M}$ that link the adjacent stations. In this case, for each adjacent station pair $m_1$ and $m_2$, two sections $(m_1, m_2)$ and $(m_2, m_1)$ exist in ${\cal E}$. Trains running in section $(m_1, m_2)$ won't influence those in the other, and vice versa. The \emph{terminal station} set ${\cal M}^\textup{ter} \subseteq {\cal M}$ specifies the stations where trains can originate and terminate. In a railway network, each station can link multiple railway sections. For each station $m\in{\cal M}$, we define a station set ${\cal M}_m = \SetCond{m'}{(m,m') \in {\cal E}}$ that includes the adjacent stations of station $m$.
	
	We are given a set $\cal T$ including the discrete timestamps within the timetable planning horizon. We are also given an original one-hour periodic timetable, which is then replicated and expanded into a multi-period original timetable to serve as input. Further, the set of lines $\cal L$ in the original timetable is given. Each \emph{line} $l \in {\cal L}$ is characterized by its OD, route, stop plan (intermediate stops)\oxc and operating frequency. Denote by ${\cal K}_l$ the set of original trains belonging to line $l \in {\cal L}$. That means, ${\cal K}_l$ contains the homogeneous trains of line $l$ distributed over all the operating hours with evenly spaced timetables. The ODs, routes, and stop plans of these original trains are fixed and cannot be modified. To maintain the regularity of train services of line $l$, it is required that at least $\left\lceil {\xi \left| {\cal K}_l \right|} \right\rceil$ trains should not be canceled, where parameter $\xi \in [0,1]$ indicates the level of periodicity of original trains. Besides, to maintain the regularity of departure and arrival times for trains of the same line, original trains' departure and arrival times at their traversed stations should not deviate from the original schedules by more than $\tau$ minutes. Further, the complete set of original trains is denoted as ${\cal K} = \bigcup\nolimits_{l \in {\cal L}} {{{\cal K}_l}}$.
	
	Extra aperiodic trains are inserted into the original timetable, and their ODs, routes, stop plans and timetables are determined within the optimization process for serving the passengers or relocating the RSUs. To limit the problem size, we are given the station paths where extra train services can possibly be operated and we are also given the departure time windows for these station paths. Multiple time windows are allowed for each station path. Generally, such time windows can be assigned to peak hours for accommodating demands, or at the beginning and ending periods to relocate the RSUs. With respect to the line plan structure, extra trains are allowed to skip at any intermediate station.
	
	To control the operating cost of the railway company, a prespecified operating cost budget $B$ is introduced to limit the total seat-kilometers of the operated original and extra train services. Besides, trains must satisfy the departure and arrival headway requirements at station boundaries to ensure operational safety, which are also referred to as track capacity requirements.
	
	\textbf{Rolling stock circulation.}
	We are given a heterogeneous fleet of RSUs, and the type set $\cal U$ specifies the types of rolling stock. An RSU of type $u \in {\cal U}$ has the seat capacity of $p_u$, and we assume that RSUs of any type $u \in {\cal U}$ are able to serve all the train services. Besides, we are given the initial number of available RSUs of type $u$ at station $m$, which is denoted by $\varpi_m^u$. A minimum transition time $\rho^\textup{conn}$ should be obeyed when an RSU finishes its previous train service and prepares to undertake the next. To consistently model the circulation of units, a larger timestamp set ${\cal T}^\textup{inv}$ is introduced that extends $\cal T$ at its ending timestamp by the transition time $\rho^\textup{conn}$. 
	
	To restrain the complexity of the problem, we only guarantee the balances of rolling stock of each type $u$ at each terminal station, i.e., we do not assign trip sequences to RSUs. This simplified version of rolling stock circulation has been adopted in previous studies \citep{Louwerse2014-MILP-Inventory,Liao2021-LR,Li2023-PeriodicExpansion}.
	
	\textbf{Passenger routing.}
	We are given a set $\cal R$ of time-dependent OD-based passenger groups, where each group $r$ consists of $g_r$ passengers with origin and destination stations $m_r^\textup{org}$ and $m_r^\textup{des}$, respectively. For each passenger group $r \in {\cal R}$, an \emph{allowable} departure time window ${\cal T}_r^\textup{all}$ and a \emph{preferred} departure time window ${\cal T}_r^\textup{pre} = \SetCond{t \in {{\cal T}_r^\textup{all}}}{t_r^\textup{min} \le t \le t_r^\textup{max}}$ are specified. Passengers can only select train services within the allowable departure time window ${\cal T}_r^\textup{all}$, and a departure shift (adaption) time cost is incurred when the departure time lies outside ${\cal T}_r^\textup{pre}$. Each group $r$ also has a latest arrival time $t_r^\textup{arr}$, by which its passengers must arrive at the destination station. 
	
	We aim to minimize total passenger travel costs, consisting of perceived travel times and penalty costs. The perceived travel time is the weighted sum of departure shift time, waiting time, in-vehicle time\oxc and transfer time, with cost coefficients $w^\textup{shift}, w^\textup{wait}, w^\textup{veh}$\oxc and $w^\textup{trans}$, respectively. The penalty cost arises when passengers in group $r \in {\cal R}$ cannot be transported due to insufficient seat capacity, in which case each unserved passenger incurs a penalty $f_r$. To ensure the level of service, we limit the number of transfers in a passenger route by $\varUpsilon$.
	
	\textbf{Hybrid periodic train timetabling problem.}
	The HPTTP integrates the above three components. The station inventory of RSUs of different types decreases and increases with the originating and terminating train services operated by the RSUs of the same type, respectively. In turn, the seat capacities of train services depend on the types of RSUs deployed. Further, each feasible combination of timetable, line plan\oxc and rolling stock circulation uniquely defines a train service network for passengers. Therefore, HPTTP aims to find such a service network design that minimizes passenger travel costs, and at the same time, ensures the consistency between passengers' routes and the service network, the consistency of rolling stock inventory and timetable, as well as operational constraints such as track capacity, operation budget\oxc and periodicity constraints.
	
	\section{Time-space network} \label{sec:ts net}
	This section introduces the time-space network representation for HPTTP, which forms the basis of model formulation. This representation has been widely adopted in train timetabling problems \citep{Brannuland1998-LR,Caprara2002-LR}, which results in integer programming models with decomposable structures that are suitable for large-scale instances.	Note that for finding pure periodic timetables which serve as inputs for our problem, the PESP model in the event-activity network \citep{Serafini-PESP} is usually used.
	
	The time-space network in this study is developed based on that of \cite{Liao2021-LR}, where train movements and inventory of RSUs are coupled. We further introduce additional station nodes that incorporate train stop decisions and transfer arcs that allow passengers to make transfers. The notations introduced in this section are also summarized in Appendix~\ref{app: tables}.

	\subsection{Flat node set}
	The \emph{flat railway network} constitutes the space dimension of the time-space network. Its node set~$\cal N$ is divided into seven subsets:
	\begin{itemize}
		\item Station stop departure nodes $d_m^{m',\textup{stop}}, m\in{\cal M}, m'\in{\cal M}_m$, that represent train departures bound for station $m'$ after dwelling at station $m$.
		\item Station stop arrival nodes $a_m^{m',\textup{stop}}, m\in{\cal M}, m'\in{\cal M}_m$, that represent train arrivals from station~$m'$ before dwelling at station $m$.
		\item Station skip departure nodes $d_m^{m',\textup{skip}}, m\in{\cal M}, m'\in{\cal M}_m$, that represent train departures bound for station~$m'$ after skipping station $m$.
		\item Station skip arrival nodes $a_m^{m',\textup{skip}}, m\in{\cal M}, m'\in{\cal M}_m$, that represent train arrivals from station~$m'$ before skipping station $m$.
		\item Station transfer nodes $n^\textup{trans}_m, m\in{\cal M}$, where passengers wait during a transfer, which can be a station platform or a waiting hall.
		\item Station inventory nodes $n^\textup{inv}_m, m\in{\cal M}^\textup{ter}$, that represent the dwell locations of RSUs. Define the inventory node set ${\cal N}^\textup{inv} = \SetCond{n^\textup{inv}_m}{m\in{\cal M}^\textup{ter}}$.
		\item Origin and destination nodes $o_r, d_r, r \in {\cal R}$, that represent the origin and destination of passenger group~$r$, respectively. Define the passenger OD node set ${\cal N}^\textup{od} = \SetCond{o_r, d_r}{r \in {\cal R}}$.
	\end{itemize}
	
	\subsection{Time-space network}
	The \emph{time-space network} is formally defined as a directed acyclic graph ${\cal G} = ({\cal V}, {\cal A})$, where $\cal V$ and $\cal A$ are the time-space vertices and arcs, respectively. The time-space paths defined in $\cal G$ thus represent the dynamic movements of trains and passengers. Incorporating the time dimension, the vertex set ${\cal V}$ is created by partially expanding $\cal N$ in time and expressed in Eq.~(\ref{eq: ts vertex}).
	\begin{equation}
	{\cal V} = {}\underbrace {{{\cal N}^\textup{od}}}_{\textrm{od nodes}}{} \cup {}\underbrace { \SetCond{(n,t)}{n \in {{\cal N}^\textup{inv}},t \in {{\cal T}^\textup{inv}}} }_{\textrm{inventory vertices}}{} \cup {}\underbrace {\SetCond{(n,t)}{n \in {\cal N} \setminus ( {{\cal N}^\textup{od}} \cup {{\cal N}^\textup{inv}}) ,t \in {\cal {\cal T}}} }_{\textrm{timetable vertices}}{}
		\label{eq: ts vertex}
	\end{equation}
	
	Further, the time-space arc set $\cal A$ is created to capture the movements of trains and passengers, which is divided into the following categories. Besides, denote by $c_a^r$ the cost of time-space arc $a$ for passenger group $r$.

	\textbf{Train section arcs.}
	Train section arcs $a \in {\cal A}^\textup{sec}$ represent train movements within railway sections. Denote the minimum running time of section $e = (m,m')$ as $\rho^\textup{sec}_e$. Due to train acceleration, when trains stop at the starting station $m$, an additional time $\rho^\textup{acc}_e$ should be added to the running time. Similarly, an additional time $\rho^\textup{dec}_e$ should be added when trains stop at station $m'$ due to deceleration. In combination, four types of section running arcs are created based on train stop plans at the two ending stations of the section, which are expressed by Eq. (\ref{eq: sec arcs}). Besides, the distance (in kilometers) of a section arc $a$ is denoted by $\beta_a$. As the route of an original train $k \in {\cal K}$ won't be changed, its running distance is denoted by $\beta_k$.
	\begin{equation}
		\begin{array}{l}
			\SetCond{(d_m^{m',\textup{skip}},a_{m'}^{m,\textup{skip}},t,t + \rho^\textup{sec}_e)}{e= (m,m') \in {\cal E},t \in {\cal T}:t + \rho^\textup{sec}_e \in {\cal T}} \\
			 {} \cup \SetCond{(d_m^{m',\textup{stop}},a_{m'}^{m,\textup{skip}},t,t + \rho^\textup{sec}_e + \rho^\textup{acc}_e)}{e = (m,m')\in {\cal E},t \in {\cal T}:t + \rho^\textup{sec}_e + \rho^\textup{acc}_e \in {\cal T}} \\
			 {} \cup \SetCond{(d_m^{m',\textup{skip}},a_{m'}^{m,\textup{stop}},t,t + \rho^\textup{sec}_e + \rho^\textup{dec}_e)}{e= (m,m') \in {\cal E},t \in {\cal T}:t + \rho^\textup{sec}_e + \rho^\textup{dec}_e \in {\cal T}} \\
			 {} \cup \SetCond{(d_m^{m',\textup{stop}},a_{m'}^{m,\textup{stop}},t,t + \rho^\textup{sec}_e + \rho^\textup{acc}_e + \rho^\textup{dec}_e)}{e= (m,m') \in {\cal E},t \in {\cal T}:t + \rho^\textup{sec}_e + \rho^\textup{acc}_e + \rho^\textup{dec}_e \in {\cal T}}
		\end{array}
		\label{eq: sec arcs}
	\end{equation}
	
	\textbf{Train dwell and passing arcs.}
	Train dwell arcs link an arrival stop vertex to a departure stop vertex in a station, and the minimum and maximum dwell times are $\underline{\rho}^\textup{dwell}_m$ and $\overline{\rho}^\textup{dwell}_m$, respectively. It is noted that no turnover occurs along the train routes; therefore, the arrival and departure nodes are not associated with the same adjacent station. The train dwell arcs are expressed in Eq. (\ref{eq : dwell arcs}).
	\begin{equation}
        \SetCond{(a_m^{m',\textup{stop}},d_m^{m'',\textup{stop}},t,t + t')}{m \in {\cal M},m',m'' \in {{\cal M}_m},m' \ne m'',t \in {\cal T},t' \in [\underline{\rho}^\textup{dwell}_m,\overline{\rho}^\textup{dwell}_m]:t + t' \in {\cal T}}
	\label{eq : dwell arcs}
	\end{equation}
	
	Train passing arcs link an arrival skip vertex to a departure skip vertex in a station. The time lengths of passing arcs are approximated as zeros. Similar to the dwell arcs, train passing arcs are expressed in Eq. (\ref{eq : passing arcs}).
	\begin{equation}
            \SetCond{(a_m^{m',\textup{skip}},d_m^{m'',\textup{skip}},t,t)}{m \in {\cal M},m',m'' \in {{\cal M}_m},m' \ne m'',t \in {\cal T}}
	\label{eq : passing arcs}
	\end{equation}

	As train section, dwell and passing arcs all indicate the in-vehicle traveling of passengers, the arc costs are calculated by multiplying coefficient $w^\textup{veh}$ by the time length of the arc.
		
	\textbf{Inventory arcs.}
	Inventory arcs represent the waiting process of RSUs at the inventory nodes at terminal station $m \in {\cal M}^\textup{ter}$. Denote the inventory node at station $m$ as $n^\textup{inv}_m$, and denote by $\sigma(t)$ the earliest timestamp after $t$, i.e., $\sigma (t) = \arg \min \SetCond{s \in {\cal T}}{s > t}$, then the inventory arcs are expressed in Eq. (\ref{eq : inventory arcs}).
	\begin{equation}
	   \SetCond{(n^\textup{inv}_m,n^\textup{inv}_m,t,\sigma (t))}{m \in {{\cal M}^\textup{ter}},t,\sigma (t) \in {\cal T}}	
		\label{eq : inventory arcs}
	\end{equation}

	\textbf{Train source arcs.}
	Train source arcs link the inventory node with a departure stop vertex at the same terminal station $m \in {\cal M}^\textup{ter}$, and represent the sign-on activities of RSUs. The source arcs are expressed in Eq. (\ref{eq : tra source arcs}).
	\begin{equation}
            \SetCond{(n^\textup{inv}_m,d_m^{m',\textup{stop}},t,t)}{m \in {{\cal M}^\textup{ter}},m' \in {{\cal M}_m},t \in {\cal T}}
		\label{eq : tra source arcs}
	\end{equation}
	
	\textbf{Train sink arcs.}
	Train sink arcs link an arrival stop vertex with the inventory node at the same terminal station $m \in {\cal M}^\textup{ter}$. The arcs represent the sign-off activities of RSUs and have the time lengths of the RSU transition time $\rho^\textup{conn}$, i.e., the RSUs won't become available for the next train service until $\rho^\textup{conn}$ minutes later. The sink arcs are expressed in Eq. (\ref{eq : tra sink arcs}).
	\begin{equation}
        \SetCond{(a_m^{m',\textup{stop}},n^\textup{inv}_m,t,t + \rho^\textup{conn} )}{m \in {{\cal M}^\textup{ter}},m' \in {{\cal M}_m},t \in {\cal T}}
	\label{eq : tra sink arcs}
	\end{equation}
	
	\textbf{Passenger transfer walking arcs.}
	Transfer walking arcs represent the walking process in passengers' transfers. The arcs link an arrival stop vertex with a transfer vertex in a station. Denote the transfer walking time at station $m$ as $\rho^\textup{trans}_m$, which is the time that passengers take to walk to the platform of the connecting train. The transfer arcs are expressed in Eq. (\ref{eq : transfer arcs}) and the arc cost is thus $w^\textup{trans}\rho^\textup{trans}_m$.
	\begin{equation}
        \SetCond{(a_m^{m',\textup{stop}},n^\textup{trans}_m,t,t + \rho^\textup{trans}_m)}{m \in {\cal M}, m' \in {\cal M}_m,t \in {\cal T}:t + \rho^\textup{trans}_m \in {\cal T}}
	\label{eq : transfer arcs}
	\end{equation}
	
	\textbf{Passenger transfer waiting arcs.}
	Transfer waiting arcs represent the passenger waiting process at the platform or the waiting hall, which happens after the passenger walking process. The transfer waiting arcs are expressed in Eq. (\ref{eq : waiting arcs}). The cost of the waiting arc is $w^\textup{wait}(\sigma(t) - t)$. It is noted that the passenger transfer time is the sum of the transfer walking time and the transfer waiting time.
	\begin{equation}
            \SetCond{(n^\textup{trans}_m,n^\textup{trans}_m,t,\sigma (t))}{m \in {\cal M},t,\sigma (t) \in {\cal T}}
		\label{eq : waiting arcs}
	\end{equation}

	\textbf{Passenger transfer boarding arcs.}
	Transfer boarding arcs represent the boarding activities of transfer passengers, which link a transfer vertex to a departure stop vertex of a station. The boarding arcs have zero time lengths and zero arc costs, and they are expressed in Eq. (\ref{eq : boarding arcs}).
	\begin{equation}
        \SetCond{(n^\textup{trans}_m, d_m^{m',\textup{stop}},t,t)}{m \in {\cal M},m' \in {\cal M}_m, t \in {\cal T}}	
	\label{eq : boarding arcs}
	\end{equation}

	\textbf{Passenger origin and destination arcs.}
	The origin arcs link the origin node with a departure stop vertex at the passenger's origin station. For consistency of expression, we use $(o_r, 0)$ to represent the static origin vertex of passenger group $r$. The origin arcs are expressed in Eq. (\ref{eq : ori arcs}). If passengers depart within the preferred time window, i.e., $t \in [t_r^\textup{min}, t_r^\textup{max}]$, the arc cost is zero. Otherwise, the cost is $w^\textup{shift}(t - t_r^\textup{max})$ if $t > t_r^\textup{max}$ or $w^\textup{shift}(t_r^\textup{min} -t)$ if $t < t_r^\textup{min}$ to account for the departure shift time.
	\begin{equation}
            \SetCond{({o_r},d_{m_r^\textup{org}}^{m',\textup{stop}},0,t)}{r \in {\cal R},m' \in {{\cal M}_{m_r^\textup{org}}},t \in {{\cal T}_r^\textup{all}}}
		\label{eq : ori arcs}
	\end{equation}
	
	Similarly, using $(d_r, 0)$ as the static destination vertex, destination arcs are expressed in Eq. (\ref{eq : des arcs}). The costs of destination arcs are zero.
	\begin{equation}
		\SetCond{(a_{m_r^\textup{des}}^{m',\textup{stop}},{d_r},t,0)}{r \in {\cal R},m' \in {{\cal M}_{m_r^\textup{des}}}, t \in {\cal T}: t \leq t_r^\textup{arr}}
	\label{eq : des arcs}
	\end{equation}

	Please note that passenger transfer nodes are introduced in our network representation in order to accurately capture passenger transfer costs. Explanations are provided in Appendix~\ref{app: transfer capture}.

	\subsection{Train and passenger time-space subnetworks.}
	After preparing the complete time-space network, this section details the subnetworks for specific trains and passengers based on the property of the problem.
	
	\textbf{Original train subnetworks.}
	For each original train $k \in {\cal K}$, the associated line information (route and intermediate stops) and its timetable are given. We construct the time-space subnetwork ${\cal G}_k = ({\cal V}_k, {\cal A}_k)$ based on the allowed schedule deviation $\tau$. For section, dwell\oxc and passing arcs, if an arc $(n,n',t,t')$ is used in the original timetable, then the arc set $\SetCond{(n,n',j,j') \in {\cal A}}{j \in {\cal T} \cap [t - \tau ,t + \tau ],j' \in {\cal T} \cap [t' - \tau ,t' + \tau ]}$ is included in ${\cal A}_k$ to ensure consistency of station paths and stop plans, while keeping deviations of departure and arrival times within $\tau$. All vertices associated with arcs in ${\cal A}_k$ are included in ${\cal V}_k$. Additionally, source and sink arcs associated with vertices in ${\cal V}_k$ are also included in ${\cal A}_k$. For convenience, denote by ${\cal A}_k^\textup{vir}$ the set of train source and sink arcs in ${\cal A}_k$.
	
	\textbf{Extra train subnetworks.}
	For extra trains, we consider only a single subnetwork ${\cal G}^\textup{ex} = ({\cal V}^\textup{ex}, {\cal A}^\textup{ex})$. For each given physical station path, time windows are given within which extra trains could depart from the terminal stations. By assuming that the train skips all intermediate stations, we can determine the earliest timestamps when it reaches specific nodes. Conversely, assuming that the train stops at all stations and dwells for the maximum dwell time, we obtain the latest timestamps. These timestamps thus define the time windows for each node, and we include arcs in ${\cal A}^\textup{ex}$ whose head and tail times both fall within these windows. The vertex set $\cal{V}^\textup{ex}$ and the set of source and sink arcs ${\cal A}^{\textup{ex}, \textup{vir}}$ are defined similarly to that of the original train subnetworks.
	
	Together, we acquire the set of train travel arcs ${{\cal A}^\textup{tr}} = \bigcup\nolimits_{k \in {\cal K}} {( {{\cal A}_k}\backslash {\cal A}_k^\textup{vir}) }  \cup ( {{\cal A}^\textup{ex}}\backslash {{\cal A}^{\textup{ex},\textup{vir}}})$ that includes all non-virtual train arcs, i.e., the active train section, dwell\oxc and passing arcs.

	\textbf{Passenger subnetworks.}
	Similar to extra trains, we can calculate for each group $r$ the earliest possible timestamps at which it can reach each node. Based on these timestamps, we incorporate relevant train travel and passenger transfer arcs whose tail times are greater than these node timestamps but smaller than the arrival time limit $t_r^\textup{arr}$ in ${\cal A}_r$. The vertex set ${\cal V}_r$, and origin and destination arcs are defined similarly to those of the trains.
	
	\subsection{Incompatible arc sets}
	To ensure the safety of operations, trains must be separated within each section. This is enforced by imposing the minimum train departure and arrival headway times between adjacent train departures and arrivals, respectively, which can also be interpreted as track capacity requirements. Generally, seven types of train headway times are implemented, which are dependent on the stop patterns of the adjacent trains. Appendix~\ref{app: headway} gives details about these headway times. We utilize a vertex occupation-based approach to model the track capacity requirements according to \cite{Yao2023-ADMM}. Denote the set of time-space vertices at the station boundary nodes as ${{\cal V}^\textup{bound}} = \SetCond{(n,t) \in {\cal V}}{n \in \SetCond{d_m^{m',\textup{skip}},d_m^{m',\textup{stop}},a_m^{m',\textup{skip}},a_m^{m',\textup{stop}}}{m \in {\cal M},m' \in {{\cal M}_m}}, t \in {\cal T}}$. Each vertex $(n,t) \in {\cal V}^\textup{bound}$ is thus treated as a resource, and all the train arcs that occupy it according to train headway requirements are included in its incompatible arc set $\phi_{n,t}$. Therefore, by requiring that no more than one arc in $\phi_{n,t}$ can be selected, the headway requirement with respect to vertex $(n,t)$ is respected. The details of the formulation of incompatible arc sets can be seen in Appendix~\ref{app: headway}.
	
	\subsection{An example}
	
\begin{figure}[htbp]
		\centering
		\FIGURE
		{\includegraphics[width=\textwidth]{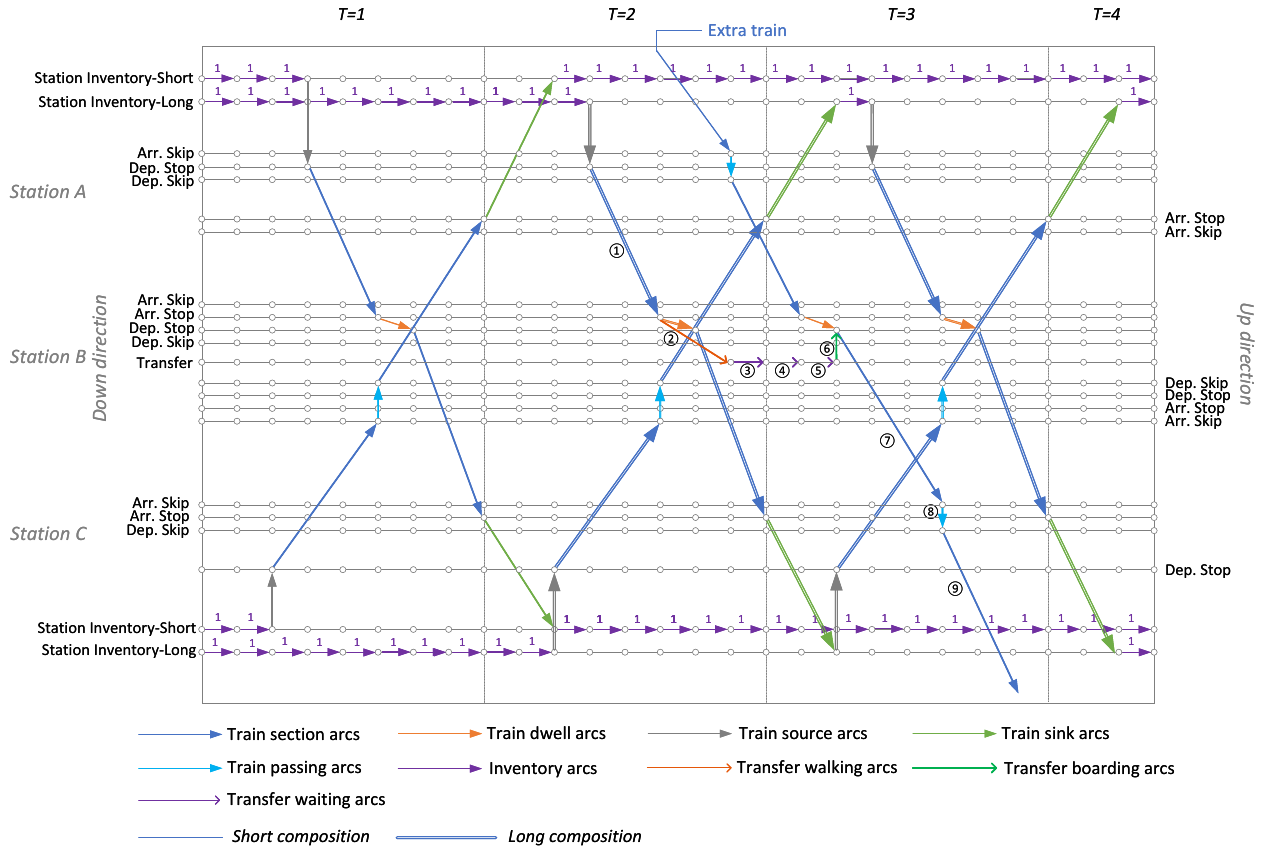}}
		{\centering{Example of the time-space network representation.}\label{fig: ts net}}
		{}
	\end{figure}

	Figure \ref{fig: ts net} gives a simple example of the time-space network representation. The example includes two types of RSUs: short-composition and long-composition ones. To improve readability, the inventory node corresponding to each type is shown separately in the figure. As vertices of some nodes are not utilized in the example, some nodes of station A and C are not shown, while station B shows the complete station nodes. Two original directed lines are operated, each containing three trains in periods 1, 2\oxc and 3. An extra train is scheduled, which originates from and terminates at locations outside the displayed network. The numbers above the inventory arcs indicate the quantities of RSUs available at the terminal stations. Besides, the circled numbers illustrate the route of a passenger from station A to station C, who first takes a train at station A during the second period, travels to station B\oxc and then makes a transfer with 5 time intervals (2 in walking and 3 in waiting) until taking the extra train for further travel. 
	
	\section{Mathematical formulation} \label{sec:model}
	The time-space network representation captures both train operations (lines and timetables) and passenger routes, enabling formulation as a dynamic service network design problem with resource constraints. Specifically, we consider two critical resources: RSUs and track resources. RSU availability is ensured through design-balanced constraints in service network design contexts \citep{Pedersen2009-DesignBal}, while track resource are managed via incompatible arc selection requirements. The subsequent section presents an arc-path formulation for the problem.
	\subsection{Arc-path formulation}
	Our arc-path formulation uses arc-based variables for train timetable and path-based variables for passenger routes. The arc-based timetabling model remains efficient due to tight time windows that limit the network scale, while the path-based passenger variables enable the use of column generation for large instances and allow the transfer requirement to be incorporated.
	
	We denote by $\delta_k^+(n,t)$ and $\delta_k^-(n,t)$ train~$k$'s arcs that direct from and into vertex $(n,t)$, respectively. Similarly, we define $\delta_\textup{ex}^+(n,t)$ and $\delta_\textup{ex}^-(n,t)$ for arcs in ${\cal A}^\textup{ex}$. Further, denote the inventory nodes at train $k$'s origin and destination station as $o_k$ and $d_k$, respectively. Denote the complete set of time-space paths of passenger group $r \in {\cal R}$ that satisfy the transfer requirement as ${\cal P}_r$, and define a binary parameter $d_p^a$ which is 1 if arc $a$ is included in path $p$ and 0 otherwise. Let $c_p^r = \sum\nolimits_{a \in {\cal A}}{d_p^a c_a^r}$ denote the travel cost of path $p$ for passenger group $r$. The following variables are considered in the model.
		\begin{itemize}
			\item[-] $x_a^{k,u} \in \{0,1\}$: equals 1 if train $k$ traverses arc $a$ using RSU of type $u$, 0 otherwise.
			\item[-] $\theta^{k,u} \in \{0,1\}$: equals 1 if train $k$ is operated with RSU of type $u$, 0 otherwise.
			\item[-] $y_a^u \in \{0,1\}$: equals 1 if an extra train traverses arc $a$ using RSU of type $u$, 0 otherwise.
			\item[-] $w_a^u \in \mathbb{Z}_{\geq 0}$: the number of RSUs of type $u$ staying at station inventory arc $a$.
			\item[-] $z_p^r \in \mathbb{R}_{\geq 0}$: the number of passengers in group $r$ that travel on path $p$.
			\item[-] $q_r \in \mathbb{R}_{\geq 0}$: the number of passengers in group $r$ that are not transported.
		\end{itemize}
	   
	The mathematical formulation for HPTTP is shown as follows:
	
	\begingroup
	\small
	\textbf{[HPTTP]}
	\begin{align}
		\min \quad& \sum\limits_{r \in {\cal R}} {\sum\limits_{p \in {{\cal P}_r}} {c_p^rz_p^r} }  + \sum\limits_{r \in {\cal R}} {{f_r}{q_r}}		\label{obj-ori-path}\\
		\textrm{s.t.} \quad & {\theta ^{k,u}} = \sum\limits_{t \in {\cal T}:({o_k},t) \in {{\cal V}_k}} {\sum\limits_{a \in \delta _k^ + ({o_k},t)} {x_a^{k,u}} } \quad \forall k \in {\cal K},u \in {\cal U}\label{con:theta x}\\
		&\sum\limits_{u \in {\cal U}} {{\theta ^{k,u}}}  \le 1 \quad \forall k \in {\cal K}
		\label{con:compo-sel}\\
		&\sum\limits_{k \in {{\cal K}_l}} {\sum\limits_{u \in {\cal U}} {{\theta ^{k,u}}} }  \ge \xi \left| {{{\cal K}_l}} \right| \quad \forall l \in {\cal L}
		\label{con:train-cancel}\\
		&\sum\limits_{a \in \delta _k^ + (n,t)} {x_a^{k,u}}  - \sum\limits_{a \in \delta _k^ - (n,t)} {x_a^{k,u}}  = 0\quad \forall k \in {\cal K},u \in {\cal U},n \in {\cal N}_k \backslash \{ o_k,d_k\},t \in {\cal T}:(n,t) \in {{\cal V}_k}
		\label{con:ori-train-balance}\\
		&\sum\limits_{a \in \delta _\textup{ex}^ + (n,t)} {y_a^u}  - \sum\limits_{a \in \delta _\textup{ex}^ - (n,t)} {y_a^u}  = 0 \quad \forall u \in {\cal U},n \in {\cal N}\backslash {{\cal N}^\textup{inv}},t \in {\cal T}:(n,t) \in {{\cal V}^\textup{ex}}\label{con:ex-train-balance}\\
		&\sum\limits_{k \in {\cal K}} {\sum\limits_{a \in {{\cal A}_k} \cap {\phi _{n,t}}} {\sum\limits_{u \in {\cal U}} {x_a^{k,u}} } }  + \sum\limits_{a \in {{\cal A}^\textup{ex}} \cap {\phi _{n,t}}} {\sum\limits_{u \in {\cal U}} {y_a^u} }  \le 1\quad \forall (n,t) \in {\cal V}^\textup{bound}\label{con:track cap}\\
		&\sum\limits_{u \in {\cal U}} {{p_u}} \left( {\sum\limits_{a \in {{\cal A}^\textup{ex}}} {{\beta_a}y_a^u}  + \sum\limits_{k \in {\cal K}} {{\beta_k}{\theta ^{k,u}}} } \right) \le B\label{con:line budget}\\
		& \sum\limits_{a \in {\delta ^ + }(n^\textup{inv}_m,0) \cap {{\cal A}^\textup{inv}}} {w_a^u}  + \sum\limits_{k \in {\cal K}} {\sum\limits_{a \in \delta _k^ + (n^\textup{inv}_m,0)} {x_a^{k,u}}  + \sum\limits_{a \in \delta _\textup{ex}^ + (n^\textup{inv}_m,0)} {y_a^u} }  = \varpi _m^u \quad \forall m \in {\cal M}^\textup{ter}, u\in{\cal U} \label{con:rsu initial}\\
		&\sum\limits_{a \in {\delta ^ - }(v) \cap {{\cal A}^\textup{inv}}} {w_a^u}  + \sum\limits_{k \in {\cal K}} {\sum\limits_{a \in \delta _k^ - (v)} {x_a^{k,u}} }  + \sum\limits_{a \in \delta _\textup{ex}^ - (v)} {y_a^u}  = \sum\limits_{a \in {\delta ^ + }(v) \cap {A^\textup{inv}}} {w_a^u}  + \sum\limits_{k \in {\cal K}} {\sum\limits_{a \in \delta _k^ + (v)} {x_a^{k,u}}  + \sum\limits_{a \in \delta _\textup{ex}^+ (v)} {y_a^u} }\nonumber\\
		&\hspace{17em} \forall v \in \SetCond{(n^\textup{inv}_m,t)}{m \in {{\cal M}^\textup{ter}},t \in {{\cal T}^\textup{inv}}\backslash \{ 0\}}, u \in {\cal U} \label{con:sta inv}\\
		&\sum\limits_{r \in {\cal R}} {\sum\limits_{p \in {{\cal P}_r}} {d_p^az_p^r} }  - \sum\limits_{k \in {\cal K}:a \in {{\cal A}_k}} {\sum\limits_{u \in {\cal U}} {{p_u}x_a^{k,u}} }  - \sum\limits_{u \in {\cal U}} {{p_u}y_a^u}  \le 0\quad \forall a \in {{\cal A}^\textup{ex}}
		\label{con:seat cap-xy-path}\\
		&\sum\limits_{r \in {\cal R}} {\sum\limits_{p \in {{\cal P}_r}} {d_p^az_p^r} }  - \sum\limits_{k \in {\cal K}:a \in {{\cal A}_k}} {\sum\limits_{u \in {\cal U}} {{p_u}x_a^{k,u}} }  \le 0\quad \forall a \in {{\cal A}^\textup{tr}} \backslash {{\cal A}^\textup{ex}}
		\label{con:seat cap-only x-path}\\		
		&\sum\limits_{p \in {{\cal P}_r}} {z_p^r + {q_r} } = {g_r} \quad \forall r\in{\cal R}
		\label{con:pax demand-path}\\
		&x_a^{k,u} \in \{0,1\} \quad \forall k \in {\cal K}, a\in{\cal A}_k, u\in{\cal U}
		\label{var x}\\
		&\theta^{k,u} \in \{0,1\} \quad \forall k \in {\cal K}, u\in{\cal U}
		\label{var theta}\\
		&y_a^u \in \{0,1\} \quad \forall a\in{\cal A}^\textup{ex}, u\in{\cal U}
		\label{var y}\\
		&w^u_a \geq 0 \; \textrm{and integer} \quad \forall u\in{\cal U}, a \in {\cal A}^\textup{inv}
		\label{var w}\\
		&z^r_p \geq 0 \quad \forall r\in{\cal R}, p \in {\cal P}_r\\
		&q_r \geq 0 \quad \forall r\in{\cal R}\label{var q}
	\end{align}
	\endgroup
	
		The objective function (\ref{obj-ori-path}) minimizes the sum of the perceived travel times and penalty costs of passengers. Constraints (\ref{con:theta x}) couple the original trains' RSU type selection variable $\theta^{k,u}$ with the source arc selection variables. Specifically, if any source arc $a$ is selected by train $k$ using RSU type $u$, then $\theta^{k,u} = 1$. Constraints (\ref{con:compo-sel}) require that at most one type of RSU can be selected by each train. If no RSUs are selected, the train is not operated. Constraints (\ref{con:train-cancel}) require that for each periodic line $l \in {\cal L}$, the proportion of the operated trains should be no smaller than the given threshold $\xi$ to guarantee the level of service periodicity.
	   
		Constraints (\ref{con:ori-train-balance}) define for each RSU type, the original train flow balance at the time-space vertices except the station inventory ones. Constraints (\ref{con:ex-train-balance}) define the extra train flow balance at time-space vertices, similarly. Constraints (\ref{con:track cap}) state the restriction on train headways. Specifically, by requiring each time-space vertex $(n,t) \in {\cal V}^\textup{bound}$ to be occupied by at most one train section arc, we can safely space all the train section arcs. Constraint (\ref{con:line budget}) specifies the operational budget limitation, which requires the seat-kilometers of operated trains be smaller than the given budget $B$. Constraints (\ref{con:rsu initial}) specify the initial RSU inventory at each terminal station at the beginning of the planning horizon, where the sum of outgoing inventory flow, original train flow and extra train flow equals the initial number of RSUs. Constraints (\ref{con:sta inv}) specify the inventory balance (design-balanced) constraints, where at each intermediate inventory vertex, the sum of incoming waiting inventory flow, original train flow\oxc and extra train flow equals that of the outgoing ones, thus guaranteeing the balance of RSUs of each type. 
		
		Constraints (\ref{con:seat cap-xy-path}) specify the seat capacity constraints for the arcs which might be traversed by extra trains. Specifically, the number of passengers on an arc should not exceed the provided seat capacity of original and extra trains. Due to headway requirements, each arc can actually be traversed by either an original or an extra train. Similarly, constraints (\ref{con:seat cap-only x-path}) specify the seat capacity constraints if the arc can only be used by original train services. Constraints (\ref{con:pax demand-path}) specify that for each passenger group, the volume of path flows plus the unrouted volume equal the number of passengers in the group. Finally, constraints (\ref{var x})--(\ref{var q}) give the domains of the variables.

		\subsection{Preprocessing techniques}
		\label{sec_prep}
		This section details some preprocessing techniques to simplify the proposed model.

		\textbf{Elimination of train arcs.}
		Due to the strict operation time window $\cal T$, all trains must complete their journeys by the operation ending time. Therefore, we compute for each train the latest feasible timestamps at each node in its subnetwork. This is achieved by performing a backward shortest path search, starting from the destination node at the ending time, and using the shortest durations of all arcs. Any arcs connected to vertices whose timestamps exceed these \emph{latest} ones are then deleted, as they cannot contribute to paths reaching the destination within the required timeframe.

		\textbf{Elimination of passenger arcs.}
		First, passengers also need to complete their travels by their latest arrival times, and the first technique can also be performed for each passenger group. Second, due to the time windows of trains, some train travel arcs will never be used. We then can delete these arcs for the passengers.

		\textbf{Upper bounding of RSU inventory variables.} The upper bound of each RSU inventory arc variable $w_a^u$ can be set as the total number of RSUs of type $u$. However, these bounds tend to be weak as they can be very large. To tighten them, one could solve $\max \SetCond{w_a^u}{\text{(\ref{con:theta x})--(\ref{var q})}}$ for all $a \in {\cal A}^\textup{inv}, u \in {\cal U}$, which however can still be difficult. Therefore, we propose a relaxation model (\ref{obj-inv-bounding})--(\ref{var x-inv-sink}). As the maximum value of $w_a^u$ is independent of passenger routing decisions, we first omit passenger-related variables and constraints. Besides, we neglect the track capacity constraints for simplification. Furthermore, for the original trains, only the source and sink arc variables are retained. Consistency between the source and sink arc selections is enforced through constraints (\ref{con:theta x}) and (\ref{con:-inv-theta x-des}). The flow balance constraints of original trains are thus dropped. The complete arc variables of extra trains are kept as they are necessary for budget calculation. The optimal objective value of the relaxation model is then used as the upper bound for variable $w_a^u$.

		\vspace{-20pt}
		\begingroup
		\small
		\begin{align} 
			\max \quad& w_a^u
			\label{obj-inv-bounding}\\
			\textrm{s.t.} \quad  &\text{(\ref{con:theta x})--(\ref{con:train-cancel})}, \text{(\ref{con:ex-train-balance})--(\ref{con:sta inv})} \nonumber \\
			&{\theta ^{k,u}} = \sum\limits_{t \in {\cal T}:({d_k},t) \in {{\cal V}_k}} {\sum\limits_{a \in \delta _k^ - ({d_k},t)} {x_a^{k,u}} } \quad \forall k \in {\cal K},u \in {\cal U}\label{con:-inv-theta x-des}\\
			&\text{(\ref{var theta})--(\ref{var w})}\nonumber\\
			&x_a^{k,u} \in \{ 0,1\} \quad \forall k \in {\cal K},a \in \SetCond{\delta _k^ + ({o_k},t)}{t \in {\cal T}:({o_k},t) \in {{\cal V}_k}},u \in {\cal U}
			\label{var x-inv-source}\\
			&x_a^{k,u} \in \{ 0,1\} \quad \forall k \in {\cal K},a \in \SetCond{\delta _k^ - ({d_k},t)}{t \in {\cal T}:({d_k},t) \in {{\cal V}_k}},u \in {\cal U}
			\label{var x-inv-sink}
		\end{align}
		\endgroup
		
		\section{Benders decomposition} \label{sec:algo}
		With a vast number of integer train arc variables and continuous passenger path variables, the arc-path model constitutes a large-scale mixed integer programming problem that becomes intractable for mathematical solvers when applied to real-world instances. For such models, particularly in the context of service network design, the Benders decomposition algorithm has been widely applied \citep{Agarwal2008-liner-scheduling,Karsten2018-BD-CG}. Therefore, we propose an efficient Benders decomposition approach that decomposes the original problem into a train master problem and a passenger subproblem, which are solved separately. By iteratively generating optimality cuts based on the solutions of the subproblem, the master problem gradually takes into account the corresponding passenger costs, leading to an effective solution procedure. 
		
		This section is organized as follows. Section \ref{sec_bd_reform} reformulates HPTTP with Benders decomposition. Section \ref{sec_cg} presents a column generation algorithm for the Benders subproblem. Section \ref{sec_imple} introduces several accelerating and cut strengthening techniques.

		\subsection{Benders reformulation}
		\label{sec_bd_reform}
		This section introduces the basic procedure of the Benders decomposition algorithm for our problem. Let ($\bm{\overline{\theta}}, \bm{\overline{w}}, \bm{\overline{x}}, \bm{\overline{y}}$) be a vector of variables that represents a feasible timetable solution satisfying constraints (\ref{con:theta x})--(\ref{con:sta inv}) and (\ref{var x})--(\ref{var w}), then given ($\bm{\overline{\theta}}, \bm{\overline{w}}, \bm{\overline{x}}, \bm{\overline{y}}$), the model reduces to the following \emph{primal Benders subproblem} (PBSP):

		\vspace{-15pt}
		\begingroup
		\small
		\begin{align}
			\textbf{[PBSP]} \quad \min \quad & \sum\limits_{r \in {\cal R}} {\sum\limits_{p \in {{\cal P}_r}} {c_p^rz_p^r} }  + \sum\limits_{r \in {\cal R}} {{f_r}{q_r}}\\
			\textrm{s.t.} \quad  &\sum\limits_{r \in {\cal R}} {\sum\limits_{p \in {{\cal P}_r}} {d_p^az_p^r} }  \le \sum\limits_{k \in {\cal K}:a \in {{\cal A}_k}} {\sum\limits_{u \in {\cal U}} {{p_u}\overline x _a^{k,u}} }  + \sum\limits_{u \in {\cal U}} {{p_u}\overline y _a^u} \quad \forall a \in {{\cal A}^\textup{ex}}
			\label{con:seat cap-PBSP-xy}\\
			&\sum\limits_{r \in {\cal R}} {\sum\limits_{p \in {{\cal P}_r}} {d_p^az_p^r} }  \le \sum\limits_{k \in {\cal K}:a \in {{\cal A}_k}} {\sum\limits_{u \in {\cal U}} {{p_u}\overline x _a^{k,u}} } \quad \forall a \in {{\cal A}^\textup{tr} \backslash {\cal A}^\textup{ex}}
			\label{con:seat cap-PBSP-only x}\\
			&\sum\limits_{p \in {{\cal P}_r}} {z_p^r + {q_r} } = {g_r} \quad \forall r\in{\cal R}
			\label{con:pax flow-PBSP}\\
			&z^r_p \geq 0 \quad \forall r\in{\cal R}, p \in {\cal P}_r
			\label{var y-PBSP}\\
			&q_r \geq 0 \quad \forall r\in{\cal R}
			\label{var z-PBSP}
		\end{align}
		\endgroup

		Let $\bm{\lambda} = (\lambda_a \mid a\in{\cal A}^\textup{tr} )$ and $\bm{\mu} = (\mu_r \mid r\in{\cal R})$ be the vectors of the dual variables associated with seat capacity constraints (\ref{con:seat cap-PBSP-xy}),(\ref{con:seat cap-PBSP-only x}) and passenger flow conservation constraints (\ref{con:pax flow-PBSP}), then the \emph{dual Benders subproblem} (DBSP) is formulated as follows:

		\vspace{-15pt}
		\begingroup
		\small
		\begin{align}
			\textbf{[DBSP]} \quad \max \quad &\sum\limits_{a \in {{\cal A}^\textup{tr}}} {\sum\limits_{k \in {\cal K}:a \in {{\cal A}_k}} {\sum\limits_{u \in {\cal U}} {{p_u}\overline x _a^{k,u}{\lambda _a}} } }  + \sum\limits_{a \in {{\cal A}^\textup{ex}}} {\sum\limits_{u \in {\cal U}} {{p_u}\overline y _a^u{\lambda _a}} }  + \sum\limits_{r \in {\cal R}} {{g_r}{\mu _r}}\label{obj con: DBSP}\\
			\textrm{s.t.} \quad & \sum\limits_{a \in {{\cal A}^\textup{tr}}} {d_p^a{\lambda _a}}  + {\mu _r} \le c_p^r \quad \forall r\in {\cal R}, p \in {\cal P}_r
			\label{con:rc-DBSP}\\
			&{\mu _r} \le {f_r} \quad \forall r\in{\cal R}
			\label{con:mu-DBSP}\\
			&\lambda_a \leq 0 \quad \forall a\in {\cal A}^\textup{tr}
			\label{con:lam-DBSP}
		\end{align}
		\endgroup

		Let $\varTheta$ denote the set of extreme points of the polyhedron defined by constraints (\ref{con:rc-DBSP})--(\ref{con:lam-DBSP}). Then, the optimal objective value of DBSP equals that of model (\ref{model:eta}), where the additional variable $\eta$ is introduced to represent the maximum value attained among all the extreme points.

		\vspace{-10pt}
		\begingroup
		\small
		\begin{equation}
			\min \; \eta \quad \textrm{s.t.} \quad \eta  \ge \sum\limits_{a \in {{\cal A}^\textup{tr}}} {\sum\limits_{k \in {\cal K}:a \in {{\cal A}_k}} {\sum\limits_{u \in {\cal U}} {{p_u}\bar x_a^{k,u}{\lambda _a}} } }  + \sum\limits_{a \in {{\cal A}^\textup{ex}}} {\sum\limits_{u \in {\cal U}} {{p_u}\bar y_a^u{\lambda _a}} }  + \sum\limits_{r \in {\cal R}} {{g_r}{\mu _r}}  \quad \forall (\bm{\lambda} ,\bm{\mu} ) \in \varTheta
			\label{model:eta}
		\end{equation}
		\endgroup
	   
		Given any feasible solution ($\bm{\overline{\theta}}, \bm{\overline{w}}, \bm{\overline{x}}, \bm{\overline{y}}$), the optimal value of $\eta$ then equals the optimal objective of PBSP due to LP duality theory. Therefore, the \emph{Benders master problem} (BMP) of the original problem is expressed in (\ref{eq_eta})--(\ref{con:opt cut}), where constraints (\ref{con:opt cut}) are called the \emph{Benders optimality cuts}.

		\vspace{-15pt}
		\begingroup
		\small
		\begin{align}
			\textbf{[BMP}(\varTheta)\textbf{]} \quad \min \quad & \eta \label{eq_eta} \\
			\text{s.t.} \quad & \text{(\ref{con:theta x})--(\ref{con:sta inv})}, \text{(\ref{var x})--(\ref{var w})} \nonumber\\
			&\eta  \ge \sum\limits_{a \in {{\cal A}^\textup{tr}}} {\sum\limits_{k \in {\cal K}:a \in {{\cal A}_k}} {\sum\limits_{u \in {\cal U}} {{p_u} x_a^{k,u}{\lambda _a}} } }  + \sum\limits_{a \in {{\cal A}^\textup{ex}}} {\sum\limits_{u \in {\cal U}} {{p_u} y_a^u{\lambda _a}} }  + \sum\limits_{r \in {\cal R}} {{g_r}{\mu _r}}  \quad \forall (\bm{\lambda} ,\bm{\mu} ) \in \varTheta
			\label{con:opt cut}
		\end{align}
		\endgroup

		It is noted that Benders feasibility cuts are not introduced, as for any given master problem solution, the subproblem always remain feasible with $q_r = g_r$ for all $r \in {\cal R}$, meaning no passengers are transported. Further, due to the potentially large number of optimality cuts (\ref{con:opt cut}), the Benders master problem solves a relaxation BMP($\overline{\varTheta}$) with a subset of optimality cuts $\overline{\varTheta} \subseteq \varTheta$. After solving BMP, DBSP is solved using the optimal solution of BMP as input. If the optimal objective value of DBSP exceeds the current $\eta$, indicating that the optimality cut derived from DBSP's solution is violated, this cut is added to BMP. The algorithm reaches optimality when $\eta$ equals DBSP's best objective, meaning all optimality cuts are satisfied.
	
		\subsection{Solving PBSP by column generation}\label{sec_cg}
		Given a solution of BMP, DBSP is then solved to generate optimality cuts. Instead of solving DBSP, we choose to solve PBSP due to its decomposable structure. PBSP is a path-based capacitated multi-commodity flow problem with a large number of path variables that cannot be enumerated a priori. Therefore, we adopt the column generation (CG) approach for PBSP where a reduced set of columns (corresponding to passengers' paths) are included in the model, and columns with negative reduced costs are added to the model to hopefully improve the solution. Specifically, after solving PBSP, one obtains the optimal dual variables $\lambda_a$ of constraints (\ref{con:seat cap-PBSP-xy}), (\ref{con:seat cap-PBSP-only x}) and $\mu_r$ of constraints (\ref{con:pax flow-PBSP}). The reduced cost $\textup{rc}_p^r$ of a path $p$ belonging to passenger group $r$ can thus be expressed by (\ref{eq: reduced cost}).

		\vspace{-10pt}
		\begingroup
		\small
		\begin{equation}
			\textup{rc}_p^r = c_p^r - \sum\limits_{a \in {{\cal A}^\textup{tr}}} {d_p^a{\lambda _a}}  - {\mu _r} = \sum\limits_{a \in {\cal A}\backslash {{\cal A}^\textup{tr}}} {d_p^ac_a^r}  + \sum\limits_{a \in {{\cal A}^\textup{tr}}} {d_p^a\left( {c_a^r - {\lambda _a}} \right)}  - {\mu _r}
			\label{eq: reduced cost}
		\end{equation}
		\endgroup
	
		To identify the promising column with the most negative reduced cost for passenger group~$r$, we should find the path $p \in {\cal P}_r$ that minimizes the reduced cost, which is equivalent to solving the shortest path problem (\ref{cg:obj})--(\ref{cg:var z}), where $\overline c _a^r = c_a^r - \lambda_a$  for $a \in {\cal A}^\textup{tr}$ and $\overline c _a^r = c_a^r$ for $a \in {\cal A} \setminus {\cal A}^\textup{tr}$. The resulting problem is also called the pricing subproblem (PP). The binary variable $z_a^r$ is introduced to represent whether arc $r$ is selected in the path. $\delta_{r}^+(n,t)$ and $\delta_{r}^-(n,t)$ denote the sets of passenger group $r$'s arcs that direct out of and into vertex $(n,t)$, respectively. Constraint (\ref{cg:transfer}) limits the number of passenger transfers in the path, where ${\cal A}^\textup{trans}_r$ represents the set of transfer walking arcs in ${\cal A}_r$. If the optimal objective is negative, then a new column derived from the optimal solution is added to the model. 

		\vspace{-13pt}
		\begingroup
		\small
		\begin{align}
			\textbf{[PP}_r\textbf{]} \quad \min \quad & {\sum _{a \in {{\cal A}_r}}}\overline c _a^rz_a^r\label{cg:obj}\\
			\textrm{s.t.} \quad & \sum\limits_{a \in \delta _r^ + (n,t)} {z_a^r}  - \sum\limits_{a \in \delta _r^ - (n,t)} {z_a^r}  = \begin{cases}
				1 &\textrm{if } (n,t) = ({o_r},0), \\
				-1 &\textrm{if } (n,t) = ({d_r},0), \\
				0 &\textrm{ otherwise}.
				\end{cases} \quad \forall (n,t)\in{\cal V}_r\\
			&\sum\limits_{a \in {\cal A}_r^\textup{trans}} {z_a^r}  \le \varUpsilon
			\label{cg:transfer}\\
			&z_a^r \in \{0,1\} \quad \forall a \in {\cal A}_r
			\label{cg:var z}
		\end{align}
		\endgroup
		To solve the shortest path problem efficiently in large-scale instances, we adopt a label correcting algorithm, where the requirement on the number of transfers is handled by treating the transfer walking arcs as resources and using approaches for resource-constrained shortest path problems \citep{Irnich2006-SPPRC}. The algorithm implemented is based on  Algorithm 2 in the work of \cite{Gudapati2022-SNDSR}.
		\subsection{Implementation strategies}
		\label{sec_imple}

		To enhance the performance of the proposed Benders decomposition algorithm, this section presents several implementation strategies, including acceleration techniques for both BMP and PBSP, along with two types of strengthened Benders optimality cuts.

		\textbf{Accelerating BMP computation.}
        As an integer program, the relaxed Benders master problem can be difficult to solve and take much computation time, which influences the algorithm's overall efficiency. It is noted that any extreme point of the subproblem's dual space yields a valid optimality cut. The solution obtained in the master problem determines only the objective of the dual subproblem and, consequently, the extreme point selected. However, it is not necessary to use an optimal solution to the (relaxed) BMP for that.
		
		The first strategy is the two-phase Benders decomposition approach of \cite{Cordeau2001-3phase}, where we firstly relax the integrality requirements of the BMP variables, and generate optimality cuts using LP optimal solutions. After a certain number of iterations, the optimality cuts in the LP phase are kept, and the integrality requirements are introduced to seek feasible solutions. 
		
		The second strategy involves not solving the integer BMP optimally but instead using suboptimal solutions with an acceptable optimality gap (AOG) to generate cuts \citep{Geoffrion1974-epsilon-Benders}. Such suboptimal solutions do not necessarily provide a lower bound for the original problem, however, during the solution of the integer program, we also obtain a lower bound that we will use in Algorithm~\ref{algo:BD}.
        To balance the efficiency and optimality of the algorithm, we propose to adaptively set the AOG of BMP. Denote the current optimality gap of the overall Benders decomposition algorithm as $\varrho$, the initial AOG of BMP as $\alpha_0$, the target AOG as $\bar{\alpha}\,(< \alpha_0)$\oxc and the decreasing rate as $\epsilon$. Then $\alpha_k$, which represents AOG of the iteration $k\,(> 1)$, is set according to Eq. (\ref{alpha set}).

		\vspace{-10pt}
		\begingroup
		\small
		\begin{equation}
			\alpha_k  = \begin{cases}
				\max \{ \bar \alpha ,{\alpha _{k - 1}} \cdot (1 - \varepsilon )\} &\textrm{if } \varrho \geq 2 \cdot \bar \alpha, \\
				1/2 \cdot \varrho &\textrm{otherwise}.
				\end{cases}
			\label{alpha set}
		\end{equation}
		\endgroup

		When the current optimality gap $\varrho$ is large, i.e., we are at the initial iterations of the algorithm, the value of $\max \{ \bar \alpha ,{\alpha _{k - 1}} \cdot (1 - \varepsilon )\}$ decreases from $\alpha_0$ to $\bar \alpha$ with the decreasing rate of $\epsilon$ per iteration. The setting of large $\alpha_k$ at initial iterations is based on the fact that BMP has too little information about the subproblem cost to be worth optimizing accurately \citep{Geoffrion1974-epsilon-Benders}. When $\alpha_k$ reaches $\bar \alpha$, the relatively small value still allows BMP to be quickly optimized and facilitates further reducing $\varrho$. When $\varrho < 2 \cdot \bar \alpha$, i.e., we are near the optimality, $\alpha_k$ is then set as $1/2 \cdot \varrho\, (<\bar \alpha)$, so as to pursue optimality. As $\varrho$ continues to decrease, the BMP is optimized more strictly. The algorithm then proceeds to optimality when $\varrho$ and $\alpha_k$ reach zero. 

		\textbf{Accelerating column generation.} Due to the vast number of passenger groups and large passenger routing subnetworks, the column generation algorithm for PBSP still requires much computation time. To accelerate computation, a hybrid column generation algorithm utilizing partial pricing and column-copying techniques is designed.

		Instead of solving all pricing subproblems in each iteration of column generation, the \emph{partial pricing technique} terminates the computation of pricing subproblems once $s < \left|\cal R \right|$ new columns are generated. These new columns are then incorporated into the restricted master problem for optimization \citep{Gamache1999-PartialPricing}. Pricing subproblems are solved in a predefined order, and the subsequent iteration resumes from the pricing subproblem where the previous iteration terminated. Although this approach may require more iterations of column generation, the overall computation time might still be reduced since each iteration now requires much less computation time.

		Second, we propose a \emph{column-copying technique} to quickly generate a large number of feasible columns in the early phase of column generation. Since the original timetable is periodic, the original passenger paths are also periodic for passengers of the same OD pair but different periods. By shifting in time an existing passenger path (column) to another period, we can obtain a feasible path for the passenger group of the same OD pair at the new period. Further, since original trains are allowed to deviate from their original schedule by a limited time interval, we can also shift the original path several minutes forward or backward to obtain more paths. During implementation, we only solve the pricing problems for passenger groups of the first period. Whenever a new column is found, we make copies of the column for each later time period and deviate time interval.

		The \emph{hybrid column generation} algorithm is then designed as follows: It first utilizes the column-copying technique, which rapidly expands the column pool. To manage the pool's size, this phase terminates when PBSP's estimated lower bound reaches 90\% of its upper bound. The lower bound is obtained by using the Lagrangian lower bound introduced in \cite{Desrosiers2005-CGPrimer}. Subsequently, using the columns generated in the first phase as input, the partial pricing technique is applied to identify promising columns until the column generation process achieves proven optimality.

		\textbf{Strengthened Benders optimality cut based on Open/Close arcs (Open-Close cut).}
		The Benders optimality cuts can be weak which results in the slow convergence of the algorithm. We hence propose a cut strengthening technique by modifying the PBSP dual variables according to \cite{Karsten2018-BD-CG}.
	   
		Given a BMP solution, let ${\cal O}^\textup{tr}$ represent the set of \emph{open} time-space arcs that have positive seat capacities provided by master variables $\bm{\overline{x}}$ and $\bm{\overline{y}}$, and ${\cal C}^\textup{tr}$ represent the set of \emph{closed} arcs with zero seat capacities. Further, let ${\cal O}^\textup{ex}$ and ${\cal C}^\textup{ex}$ represent sets of open and closed arcs in ${\cal A}^\textup{ex}$, respectively. Using the new notations, the objective function (\ref{obj con: DBSP}) of DBSP can be expressed by Eq. (\ref{DBSP obj: open close})
		\begingroup
		\small
		\begin{equation}
			\max \sum\limits_{u \in {\cal U}} {{p_u}\left\{ {\sum _{a \in {{\cal O}^\textup{tr}}\backslash {{\cal O}^\textup{ex}}}}{\lambda _a}\sum\limits_{k \in {\cal K}:a \in {{\cal A}_k}} {\bar x_a^{k,u}}  + \sum\limits_{a \in {{\cal O}^\textup{ex}}} {{\lambda _a}\left( {\sum\limits_{k \in {\cal K}:a \in {{\cal A}_k}} {\bar x_a^{k,u}}  + \bar y_a^u} \right)}  + {\sum _{a \in {{\cal C}^\textup{tr}}\backslash {{\cal C}^\textup{ex}}}}0{\lambda _a} + \sum\limits_{a \in {{\cal C}^\textup{ex}}} {0{\lambda _a}} \right\}}
			+ \sum\limits_{r \in {\cal R}} {{g_r}{\mu _r}}
			\label{DBSP obj: open close}
		\end{equation}
		\endgroup
	
		It is noted that the values of $\lambda_a$, $a\in {\cal C}^\textup{tr}$, do not affect the objective function. Denote $\lambda_a^*$ and $\mu_r^*$ as the optimal solution of DBSP, then as long as $\lambda_a^*$, $a \in {\cal O}^\textup{tr}$, and $\mu_r^*$, $r\in {\cal R}$, are fixed, there can be many optimal solutions by using different $\lambda_a$ for $a\in {\cal C}^\textup{tr}$ that satisfy constraints (\ref{con:rc-DBSP}) and (\ref{con:lam-DBSP}). Therefore, to tighten the optimality cut, we would like $\lambda_a$, $a\in {\cal C}^\textup{tr}$, to be as large as possible, which can be realized by solving the following cut strengthening problem:

		\vspace{-25pt}
		\begingroup
		\small
		\begin{align}
			\max \quad & \sum\limits_{a \in {{\cal C}^\textup{tr}}} {{\lambda _a}}\\
			\textrm{s.t.} \quad & \sum\limits_{a \in {{\cal C}^\textup{tr}}} {d_p^a{\lambda _a}}  \le c_p^r - \sum\limits_{a \in {{\cal O}^\textup{tr}}} {d_p^a\lambda _a^*}  - \mu _r^*\quad \forall r \in {\cal R},p \in {{\cal P}_r}
			\label{con:rc-DBSP-strengthen}\\
			&\lambda_a \leq 0 \quad \forall a\in {\cal C}^\textup{tr}
			\label{con:lam-DBSP-strengthen}
		\end{align}
		\endgroup

		Then, the values of $\lambda_a$, $a\in {\cal C}^\textup{tr}$ produced by the cut strengthening problem are thus combined with $\lambda_a^*$, $a \in {\cal O}^\textup{tr}$, and $\mu_r^*$, $r\in {\cal R}$ to derive the strengthened optimality cut. To consistently utilize our column generation framework, we also solve the dual of the cut strengthening problem similar to \cite{Karsten2018-BD-CG}, and details can be seen in Appendix~\ref{app: oc cut}.

		\textbf{Pareto-optimal cut.} In addition to the Open-Close cut, we further propose another strengthened optimality cut. \cite{Magnanti1981-pareto} proposed a method to generate a cut that is not dominated by other cuts, which is defined as the ``Pareto-optimal cut''. Such strengthened cut could improve the convergence of the algorithm and reduce the number of Benders iterations. To find such cut, one needs to solve an auxiliary problem using a core point which is within the relative interior of the convex hull of all BMP solutions. Define ($\hat{\bm{x}}, \hat{\bm{y}}$) as the core point and $Q({\bm{\lambda} ^*},{\bm{\mu}^*})$ as the optimal objective value of the DBSP defined at the current BMP solution ($\bar{\bm{x}}, \bar{\bm{y}}$). We can solve the following model (\ref{obj con: DBSP-pareto})--(\ref{con: pareto}) to find the Pareto-optimal cut.
		
		\vspace{-20pt}
		\begingroup		\small
		\begin{align}
			\max \quad &\sum\limits_{a \in {{\cal A}^\textup{tr}}} {\sum\limits_{k \in {\cal K}:a \in {{\cal A}_k}} {\sum\limits_{u \in {\cal U}} {{p_u}\hat x _a^{k,u}{\lambda _a}} } }  + \sum\limits_{a \in {{\cal A}^\textup{ex}}} {\sum\limits_{u \in {\cal U}} {{p_u}\hat y _a^u{\lambda _a}} }  + \sum\limits_{r \in {\cal R}} {{g_r}{\mu _r}}\label{obj con: DBSP-pareto}\\
			\textrm{s.t.} \quad & \sum\limits_{a \in {{\cal A}^\textup{tr}}} {d_p^a{\lambda _a}}  + {\mu _r} \le c_p^r \quad \forall r\in {\cal R}, p \in {\cal P}_r
			\label{con:rc-DBSP-pareto}\\
			&{\mu _r} \le {f_r} \quad \forall r\in{\cal R}
			\label{con:mu-DBSP-pareto}\\
			&\lambda_a \leq 0 \quad \forall a\in {\cal A}^\textup{tr}
			\label{con:lam-DBSP-pareto}\\
			& \sum\limits_{a \in {{\cal A}^\textup{tr}}} {\sum\limits_{k \in {\cal K}:a \in {{\cal A}_k}} {\sum\limits_{u \in {\cal U}} {{p_u}\bar x_a^{k,u}{\lambda _a}} } }  + \sum\limits_{a \in {{\cal A}^\textup{ex}}} {\sum\limits_{u \in {\cal U}} {{p_u}\bar y_a^u{\lambda _a}} }  + \sum\limits_{r \in {\cal R}} {{g_r}{\mu _r}}  = Q({\bm{\lambda} ^*},{\bm{\mu}^*})
			\label{con: pareto}
		\end{align}
		\endgroup

		Constraint (\ref{con: pareto}) ensures that the solution to the auxiliary problem is out of the optimal solutions to the DBSP. Then, objective function (\ref{obj con: DBSP-pareto}) seeks to find the solution that provides the maximum objective value at the core point, and such solution is thus used to generate the Pareto-optimal cut. As the core point is nontrivial to obtain, we use the following core point update mechanism to approximate the core point according to \cite{WangYD2019-liner}:

		\vspace{-10pt}
		\begingroup
		\small
		\begin{equation}
			(\hat{\bm{x}}, \hat{\bm{y}})_i = 
			\begin{cases}
			(\bar{\bm{x}}, \bar{\bm{y}}), & i = 1, \\
			\left[ (i-1) \cdot (\hat{\bm{x}}, \hat{\bm{y}})_{i-1} + (\bar{\bm{x}}, \bar{\bm{y}}) \right] / i, & i > 1, 
			\end{cases}
		\end{equation}
		\endgroup

		where $(\bar{\bm{x}}, \bar{\bm{y}})$ is the current BMP solution, $i$ is the iteration index\oxc and $(\hat{\bm{x}}, \hat{\bm{y}})_{i}$ is the core point at $i$th iteration. At the first iteration, the core point is approximated as the first obtained solution in the BMP, and then approximated to the middle of all such solutions obtained in the following iterations.

		Still, the auxiliary problem contains an increasingly large number of constraints (\ref{con:rc-DBSP-pareto}). To consistently utilize our column generation algorithm, we propose to solve the dual of the auxiliary problem. Denote by $z_p^r$, $q_r$, and $\kappa$ the dual variables of constraints (\ref{con:rc-DBSP-pareto}), (\ref{con:mu-DBSP-pareto}), and (\ref{con: pareto}), respectively, then the dual problem is formulated as (\ref{obj:PBSP-pareto})--(\ref{var kappa-pareto}). By starting from a limited set of paths $\cal{P}_r$, $r \in \cal{R}$, new paths are generated through solving the pricing subproblems introduced in Section \ref{sec_cg}.

		\vspace{-20pt}
		\begingroup
		\small
		\begin{align}
			\min \quad & \sum\limits_{r \in {\cal R}} {\sum\limits_{p \in {{\cal P}_r}} {c_p^rz_p^r} }  + \sum\limits_{r \in {\cal R}} {{f_r}{q_r}} + Q({\bm{\lambda} ^*},{\bm{\mu}^*})\kappa \label{obj:PBSP-pareto}\\
			\textrm{s.t.} \quad  &\sum\limits_{r \in {\cal R}} {\sum\limits_{p \in {{\cal P}_r}} {d_p^az_p^r} } + \left( \sum\limits_{k \in {\cal K}:a \in {{\cal A}_k}} {\sum\limits_{u \in {\cal U}} {{p_u}\overline x _a^{k,u}} }  + \sum\limits_{u \in {\cal U}} {{p_u}\overline y _a^u}\right)\kappa \le \sum\limits_{k \in {\cal K}:a \in {{\cal A}_k}} {\sum\limits_{u \in {\cal U}} {{p_u}\hat x _a^{k,u}} }  + \sum\limits_{u \in {\cal U}} {{p_u}\hat y _a^u} \quad \forall a \in {{\cal A}^\textup{ex}}
			\label{con:seat cap-PBSP-xy-pareto}\\
			&\sum\limits_{r \in {\cal R}} {\sum\limits_{p \in {{\cal P}_r}} {d_p^az_p^r} } + \sum\limits_{k \in {\cal K}:a \in {{\cal A}_k}} {\sum\limits_{u \in {\cal U}} {{p_u}\overline x _a^{k,u}}}\kappa \le \sum\limits_{k \in {\cal K}:a \in {{\cal A}_k}} {\sum\limits_{u \in {\cal U}} {{p_u}\hat x _a^{k,u}} } \quad \forall a \in {{\cal A}^\textup{tr} \backslash {\cal A}^\textup{ex}}
			\label{con:seat cap-PBSP-only x-pareto}\\
			&\sum\limits_{p \in {{\cal P}_r}} {z_p^r} + {q_r} + {g_r}\kappa = {g_r} \quad \forall r\in{\cal R}
			\label{con:pax flow-PBSP-pareto}\\
			&z^r_p \geq 0 \quad \forall r\in{\cal R}, p \in {\cal P}_r
			\label{var y-PBSP-pareto}\\
			&q_r \geq 0 \quad \forall r\in{\cal R}
			\label{var z-PBSP-pareto}\\
			&\kappa \in \mathbb{R} 
			\label{var kappa-pareto}
		\end{align}
		\endgroup

	Algorithm \ref{algo:BD} outlines the proposed Benders decomposition algorithm.

	\begin{algorithm} 
		\scriptsize
		\caption{Benders decomposition for HPTTP} \label{algo:BD}
		
		\SetKwInOut{Input}{Input} 
		\SetKwInOut{Initialize}{Initialize} 
		
		\Input{Time-space network, LP maximum iteration number $K_\textup{LP}$, time limit $T$, acceptable optimality gap $\textup{Gap}$;}
		\Initialize{Current iteration number $k = 0$, global lower and upper bounds $\textup{LB}=0, \textup{UB}=+\infty$, current optimality gap $\varrho = 100\%$}
		
		Reformulate the model using Benders decomposition \hfill \tcp{Section \ref{sec_bd_reform}}
		Relax integrality requirements of BMP and obtain LP-BMP \hfill  \tcp{Two-phase Benders decomposition}
		
		\While{$k < K_\textup{LP}$ and $T$ not reached\hfill \tcp{LP phase (not necessarily optimal)}}{ 
			Solve LP-BMP to optimality\;
			Solve PBSP and generate optimality cuts; \hfill \tcp{Section \ref{sec_cg}}
			Generate strengthened optimality cuts; \hfill  \tcp{Section \ref{sec_imple}}
			Add optimality cuts to LP-BMP, $k = k+1$\;
		}
		
		Re-introduce integrality requirements and obtain MIP-BMP, set iteration number $k=0$\;
		Add the optimality cuts of LP-BMP to MIP-BMP\;
		
		\While{$\varrho\ge \textup{Gap}$ and $T$ not reached\hfill \tcp{IP phase}}{ 
			Solve MIP-BMP to the acceptable optimality gap, let $\textup{LB}^k$ be its lower bound returned from the solver\;
			Solve PBSP and generate optimality cuts, let $\textup{UB}^k$ be the objective of PBSP; \hfill \tcp{Section \ref{sec_cg}}
			Generate strengthened optimality cuts; \hfill \tcp{Section \ref{sec_imple}}
			Add optimality cuts to MIP-BMP\;
			$\textup{UB} = \min\{\textup{UB}^k, \textup{UB}\}$, $\textup{LB} = \max\{\textup{LB}^k, \textup{LB}\}$, $\varrho = (\textup{UB} - \textup{LB}) / \textup{UB}$\;
		}
		\end{algorithm}
	\section{Routing a subset of passenger groups}\label{sec: subset routing}
	It is observed that the large number of passenger groups and the resulting extensive set of passenger variables pose significant challenges in solving large-scale instances. Drawing inspiration from \cite{Schiewe2020-PR} who fixes a subset of passengers along their shortest paths and only routes the rest during optimization, we propose a novel model tailored to our problem setting, solved using a customized Benders decomposition algorithm.
	
	To implement such a technique in our problem setting, two differences should be taken care of. First, in the periodic setting, line-based transfers of the fixed passenger routes can always be guaranteed in the optimization phase since the connecting line is regularly operated. However, in our problem setting with an acyclic time horizon, the feasibility of the train-based passenger transfers still depends on the consistency between the chosen arcs of the associated pair of trains, which requires introducing complex logical constraints. Therefore, we choose to only fix direct passengers to avoid this impact. Second, as the seat capacity is limited in our study, we need to take into account the occupied seats of the fixed passengers, and also the possibility of letting fixed passengers to be unserved when the capacity is limited, so as to increase the flexibility of the problem.

	\subsection{Model formulation}
	In this section we formulate the problem HPTTP-Passenger Subset Routing (PSR). Denote by $\cal {R}^\textup{F}$ the set of passenger groups whose routes are fixed and by $\cal{R}^\textup{F}_k$ the set of fixed passenger groups that are assigned to the original train $k$. Let $k_r$ represent the original train to which the fixed passenger group $r$ is assigned. We introduce a binary variable $\chi_r$ for $r \in \cal{R}^\textup{F}$, which equals 1 if fixed passenger group~$r$ is served and 0 otherwise. Due to the time window requirements of fixed passenger group $r \in {\cal R}$, only when train $k_r$ departs from group $r$'s origin station within its departure time window, and also arrives at its destination station before the latest arrival time, then group~$r$ can be served. In other words, only when eligible time-space arcs of the group's first and last sections are chosen by train $k_r$, then the group can be served. Therefore, for each fixed passenger group $r$, we denote by $\cal{A}^\textup{ori}_r \subseteq \cal{A}_r$ and $\cal{A}^\textup{des}_r\subseteq \cal{A}_r$ the sets of qualified section arcs of the passenger group's first and last section, respectively.

	HPTTP-PSR is formulated as a mixed-integer quadratically constrained quadratic programming model shown as follows:
	
	\begingroup
	\small
	\textbf{[HPTTP-PSR]}
	\begin{align}
		\min \quad& \sum\limits_{k \in {\cal K}} {\sum\limits_{u \in {\cal U}} {\sum\limits_{a \in {{\cal A}_k}} {\sum\limits_{r \in {\cal R}_k^\textup{F}:a \in {{\cal A}_r}} {c_a^r{g_r}{\chi _r}x_a^{k,u}} } } }  + \sum\limits_{r \in {{\cal R}^\textup{F}}} {{f_r}{g_r}\left( {1 - {\chi _r}} \right)} + \sum\limits_{r \in {\cal R}\setminus {\cal R}^\textup{F}} {\sum\limits_{p \in {{\cal P}_r}} {c_p^rz_p^r} }  + \sum\limits_{r \in {\cal R}\setminus {\cal R}^\textup{F}} {{f_r}{q_r}}  \label{obj-fix-path}\\
		\textrm{s.t.} \quad & \text{(\ref{con:theta x})--(\ref{con:sta inv})}, \text{(\ref{var x})--(\ref{var w})} \nonumber \\
		& {\chi _r} \le \sum\limits_{a \in {\cal A}_r^\textup{ori} \cap {{\cal A}_{{k_r}}}} {\sum\limits_{u \in {\cal U}} {x_a^{{k_r},u}} } \quad \forall r \in {\cal {\cal R}}^{F} \label{psr-ori-consistent}\\
		& {\chi _r} \le \sum\limits_{a \in {\cal A}_r^\textup{des} \cap {{\cal A}_{{k_r}}}} {\sum\limits_{u \in {\cal U}} {x_a^{{k_r},u}} } \quad \forall r \in {\cal {\cal R}}^{F} \label{psr-des-consistent} \\
		&\sum\limits_{r \in {\cal R}\backslash {{\cal R}^\textup{F}}} {\sum\limits_{p \in {{\cal P}_r}} {d_p^az_p^r} }  - \sum\limits_{k \in {\cal K}:a \in {{\cal A}_k}} {\sum\limits_{u \in {\cal U}} {\left( {{p_u} - \sum\limits_{r \in {{\cal R}^\textup{F}_k}:a\in{\cal A}_r} {g_r{\chi _r}} } \right)x_a^{k,u}} }  - \sum\limits_{u \in {\cal U}} {{p_u}y_a^u}   \le 0\quad \forall a \in {{\cal A}^\textup{ex}}
		\label{con:seat cap-xy-path-fix}\\
		&\sum\limits_{r \in {\cal R}\backslash {{\cal R}^\textup{F}}} {\sum\limits_{p \in {{\cal P}_r}} {d_p^az_p^r} }  - \sum\limits_{k \in {\cal K}:a \in {{\cal A}_k}} {\sum\limits_{u \in {\cal U}} {\left( {{p_u} - \sum\limits_{r \in {{\cal R}^\textup{F}_k}:a\in{\cal A}_r} {g_r{\chi _r}} } \right)x_a^{k,u}} }  \le 0\quad \forall a \in {{\cal A}^\textup{tr}} \backslash {{\cal A}^\textup{ex}}
		\label{con:seat cap-only x-path-fix}\\		
		&\sum\limits_{p \in {{\cal P}_r}} {z_p^r + {q_r} } = {g_r} \quad \forall r\in{\cal R}\setminus {\cal R}^\textup{F}
		\label{con:pax demand-path-fix}\\
		&\chi_r \in \{0,1\} \quad \forall r \in {\cal R}^\textup{F} \label{var chi} \\
		&z^r_p \geq 0 \quad \forall r\in{\cal R}\backslash {{\cal R}^\textup{F}}, p \in {\cal P}_r\label{var z-pbsp-fix}\\
		&q_r \geq 0 \quad \forall r\in{\cal R}\backslash {{\cal R}^\textup{F}}\label{var q-pbsp-fix}
	\end{align}
	\endgroup

	In the objective function (\ref{obj-fix-path}), the first term denotes the fixed passengers' routing costs, which are now represented by the arc costs of trains. Specifically, the cost of each original train arc $a \in {\cal A}_k$ is the sum of the arc costs of the fixed passengers assigned to this arc. For group $r \in {\cal R}_k^\textup{F}$, we gather the group's cost on the arc if the arc can be used by the group ($a \in{\cal A}_r$) and also the group is served ($\chi_r = 1$). The second term is thus the penalty cost of passenger group $r \in \cal{R}^\textup{F}$ if it is unserved. The last two terms represent the routing and penalty costs of unfixed passenger groups. 

	Constraints (\ref{psr-ori-consistent}) and (\ref{psr-des-consistent}) impose the consistency between the routing decision of fixed passengers and the arc selections of their assigned trains. Specifically, only when the departing and ending section arcs of passenger group $r \in {\cal R}^\textup{F}$ are both selected by its assigned train $k_r$, then the group can be served. If train $k_r$ is not operated, i.e., none of its arcs are selected at all, then also the group cannot be served. Constraints (\ref{con:seat cap-xy-path-fix}) and (\ref{con:seat cap-only x-path-fix}) specify the seat capacity requirements, where train $k \in {\cal K}$'s capacity is now influenced by the routing decision of fixed passengers. For each fixed passenger group $r \in {\cal R}_k^\textup{F}$, we deduct train $k$'s capacity by the group's volume $g_r$ on arc $a \in {\cal A}_k$ if the arc can be used by the group $a$ ($a \in {\cal A}_r$) and also the group is served ($\chi_r = 1$).

	\subsection{Benders decomposition for HPTTP-PSR}
	As HPTTP-PSR is still difficult to solve in large instances, we choose to consistently utilize our Benders decomposition approach. Let the bold letters represent the vectors of variables, it is noted that if $\bm{\chi}, \bm{x}, \bm{y}$ are fixed as $\bm{\bar{\chi}}, \bm{\bar x}, \bm{\bar y}$, then we can still get the linear primal Benders subproblem PBSP-PSR (\ref{obj:pbsp-psr})--(\ref{con:seat cap-PBSP-only x-PSR}), which is suitable for Benders decomposition framework. Therefore, the decisions of train arc selection and the routing of fixed passengers are now addressed in the master problem, while the routing of unfixed passengers are addressed in the subproblem.

	\begingroup
	\small
	\textbf{[PBSP-PSR]}
	\begin{align}
		\min \quad & \sum\limits_{r \in {\cal R}\setminus {\cal R}^\textup{F}} {\sum\limits_{p \in {{\cal P}_r}} {c_p^rz_p^r} }  + \sum\limits_{r \in {\cal R}\setminus {\cal R}^\textup{F}} {{f_r}{q_r}}\label{obj:pbsp-psr}\\
		\textrm{s.t.} \quad  &(\ref{con:pax demand-path-fix}), (\ref{var z-pbsp-fix}), (\ref{var q-pbsp-fix}) \nonumber  \\
		&\sum\limits_{r \in {\cal R}\setminus {\cal R}^\textup{F}} {\sum\limits_{p \in {{\cal P}_r}} {d_p^az_p^r} }  \le \sum\limits_{k \in {\cal K}:a \in {{\cal A}_k}} {\sum\limits_{u \in {\cal U}} {\left( {{p_u} - \sum\limits_{r \in {{\cal R}^\textup{F}_k:a\in{\cal A}_r}} {g_r{\bar \chi _r}} } \right)\bar x_a^{k,u}} }  + \sum\limits_{u \in {\cal U}} {{p_u}\bar y_a^u} \quad \forall a \in {{\cal A}^\textup{ex}}
		\label{con:seat cap-PBSP-xy-PSR}\\
		&\sum\limits_{r \in {\cal R}\setminus {\cal R}^\textup{F}} {\sum\limits_{p \in {{\cal P}_r}} {d_p^az_p^r} }  \le \sum\limits_{k \in {\cal K}:a \in {{\cal A}_k}} {\sum\limits_{u \in {\cal U}} {\left( {{p_u} - \sum\limits_{r \in {{\cal R}^\textup{F}_k:a\in{\cal A}_r}} {g_r{\bar \chi _r}} } \right)\bar x_a^{k,u}} }  \quad \forall a \in {{\cal A}^\textup{tr} \backslash {\cal A}^\textup{ex}}
		\label{con:seat cap-PBSP-only x-PSR}
	\end{align}
	\endgroup

	Based on the previous introduction of Benders optimality cuts in Section \ref{sec_bd_reform}, we can get the Benders master problem BMP-PSR (\ref{obj:bmp-psr})--(\ref{con:opt cut-psr}). It is noted that quadratic terms both exist in the objective function and optimality cuts, while they are easy to be linearized for the integer Benders master problem and therefore the details are not shown here. 

	\begingroup
	\small
	\textbf{[BMP-PSR]}
	\begin{align}
		\min \quad &\sum\limits_{k \in {\cal K}} {\sum\limits_{u \in {\cal U}} {\sum\limits_{a \in {{\cal A}_k}} {\sum\limits_{r \in {\cal R}_k^\textup{F}:a \in {{\cal A}_r}} {c_a^r{g_r}{\chi _r}x_a^{k,u}} } } }  + \sum\limits_{r \in {{\cal R}^\textup{F}}} {{f_r}{g_r}\left( {1 - {\chi _r}} \right)} + \eta \label{obj:bmp-psr} \\
		\text{s.t.} \quad & \text{(\ref{con:theta x})--(\ref{con:sta inv})}, \text{(\ref{var x})--(\ref{var w})}, (\ref{psr-ori-consistent}), (\ref{psr-des-consistent}), (\ref{var chi}) \nonumber\\
		&\eta  \ge \sum\limits_{a \in {{\cal A}^\textup{tr}}} {\sum\limits_{k \in {\cal K}:a \in {{\cal A}_k}} {\sum\limits_{u \in {\cal U}} {\left( {{p_u} - \sum\limits_{r \in {{\cal R}^\textup{F}_k}:a\in{\cal A}_r} {g_r{\chi _r}} } \right)x_a^{k,u}{\lambda _a}} } }  + \sum\limits_{a \in {{\cal A}^\textup{ex}}} {\sum\limits_{u \in {\cal U}} {{p_u} y_a^u{\lambda _a}} }  + \sum\limits_{r \in {\cal R}} {{g_r}{\mu _r}} \quad \forall (\bm{\lambda} ,\bm{\mu} ) \in \varTheta 		\label{con:opt cut-psr}
	\end{align}
	\endgroup

	It is noted that in the resulting timetable, some fixed passengers may have more favorable routes than  those initially assigned. Therefore, a post-optimization step is conducted, in which the obtained timetable is fixed, while vehicle circulation and passenger routing for \emph{all} passenger groups are jointly re-optimized. This step is implemented by re-solving the HPTTP with only the selected time-space arcs in the obtained timetable and without the extra trains.

	\section{Numerical results} \label{sec:numerical}
	
	In this section, we validate the proposed Benders decomposition algorithm. The algorithm is coded in C++ on a personal computer with 3.2GHz CPU and 64GB RAM, and Gurobi 11.0 is used as the LP and MIP solver. The integer BMP is directly solved by Gurobi, and we retrieve the best bound value from Gurobi when BMP is not solved to optimality.
	
	\subsection{Small case}
	The small case is designed based on the \texttt{Toy} dataset in the LinTim project \citep{LinTim} and the details of the case description can be seen in Appendix~\ref{app: case toy}. We generate instances using different problem settings. Two levels of operating budgets are considered, where b1 and b2 represent the budget of 400 and 500 thousand seat-kilometers, respectively. Three levels of original timetable deviations are considered, where o1, o2\oxc and o3 represent the deviation time windows of 0, 4\oxc and 8 minutes, respectively. Two time windows for extra trains are considered, where e1 and e2 represent the deviation of 0 and 4 minutes, respectively. Two periodicity levels p1 and p2 are considered, which represent the values 1 and 0.6 of $\xi$, respectively. Further, ``-t'' represents whether the passenger transfer requirement is enforced. If so, the maximum number of transfers is set as one. Together, 22 instances are generated as shown in the ``Instance'' column in Table \ref{tab: small UB LB}.
	
	\subsubsection{Benefits of the two-phase approach and strengthened optimality cuts.}
    In this experiment, we implement the algorithm with six computation settings, where ``Standard'' represents the traditional Benders decomposition and ``LP'' represents the two-phase approach. In ``LP-Pareto'', the two-phase Benders decomposition algorithm using Pareto-optimal cuts is performed, and in ``LP-Pareto-2'' both the Pareto-optimal cut and the original Benders optimality cut are added in each iteration.
    Finally, ``LP-OC'' and ``LP-OC-2'' are similarly defined with respect to the Open-Close cut.
    To expedite computation, the LP phase terminates when the optimality gap is within 1 \%. By experiments, $\alpha_0, \bar \alpha$\oxc and $\epsilon$ are set as 4\%, 1\%\oxc and 5\%, respectively. The computation time limit is set as 900s.
	
	The upper bounds (UB) and optimality gaps (Gap(\%)) are reported in Table \ref{tab: small UB LB}. If the optimality is reached, the computation time in seconds (CPU(s)) is reported instead. The bold numbers indicate the best upper bounds and optimality gaps in each instance. 
	\begin{table}[htb] \tiny
		\centering
		\caption{Experiment result with different implementation strategies of Benders decomposition}
		\begin{tabular}{ccccccccccccc}
			\toprule
			\multirow{3}[0]{*}{Instance} & \multicolumn{2}{c}{Standard} & \multicolumn{2}{c}{LP} & \multicolumn{2}{c}{LP-Pareto} & \multicolumn{2}{c}{LP-Pareto-2} & \multicolumn{2}{c}{LP-OC} & \multicolumn{2}{c}{LP-OC-2} \\
			\cmidrule(lr){2-3} \cmidrule(lr){4-5}\cmidrule(lr){6-7}\cmidrule(lr){8-9}\cmidrule(lr){10-11}\cmidrule(lr){12-13}
				  & \multirow{2}[0]{*}{UB} & CPU(s) & \multirow{2}[0]{*}{UB} & CPU(s) & \multirow{2}[0]{*}{UB} & CPU(s) & \multirow{2}[0]{*}{UB} & CPU(s) & \multirow{2}[0]{*}{UB} & CPU(s) & \multirow{2}[0]{*}{UB} & CPU(s) \\
				  &       & /Gap(\%) &       & /Gap(\%) &       & /Gap(\%) &       & /Gap(\%) &       & /Gap(\%) &       & /Gap(\%) \\
			\midrule
			o1p1b1e1 & \textbf{62007} & \textbf{6s} & 62007 & 7s    & 62007 & 12s   & 62007 & 11s   & 62007 & 14s   & 62007 & 16s \\
			o1p1b2e1 & \textbf{61861} & \textbf{3s} & 61861 & 4s    & 61861 & 8s    & 61861 & 10s   & 61861 & 9s    & 61861 & 7s \\
			o1p2b1e1 & 61114 & 2.78  & 61099 & 1.67  & 61149 & 2.08  & 61149 & 2.12  & 61149 & 2.12  & \textbf{61096} & \textbf{1.23 } \\
			o1p2b2e1 & \textbf{60666} & \textbf{45s} & 60666 & 47s   & 60666 & 76s   & 60666 & 85s   & 60666 & 76s   & 60666 & 74s \\
			o2p1b1e2 & 59663 & 2.08  & 59669 & 1.05  & 59645 & 0.98  & \textbf{59603} & \textbf{0.63 } & 59678 & 1.13  & 59844 & 1.04  \\
			o2p1b2e2 & 59027 & 0.02  & 59027 & 0.02  & \textbf{59027} & \textbf{168s} & 59051 & 0.11  & 59027 & 0.24  & 59027 & 354s \\
			o2p1b1e1 & 59727 & 0.97  & 59751 & 0.28  & 59727 & 0.11  & \textbf{59727} & \textbf{0.09 } & 59754 & 0.72  & 59799 & 0.84  \\
			o2p1b2e1 & \textbf{59364} & \textbf{89s} & 59364 & 232s  & 59364 & 184s  & 59364 & 169s  & 59364 & 0.02  & 59364 & 154s \\
			o2p2b1e1 & 61350.5 & 24.51  & \textbf{58650} & 3.98  & 58912 & 4.88  & 58793.5 & \textbf{3.52 } & 59210 & 5.55  & 58928 & 5.09  \\
			o2p2b2e1 & 58616 & 6.74  & 58299.5 & 2.39  & \textbf{58275.5} & \textbf{1.78 } & 58323.5 & 2.51  & 58434.5 & 3.16  & 58415.5 & 2.92  \\
			o2p2b1e2 & 61427 & 36.04  & 59407.5 & 6.67  & 58791 & 5.68  & 58586.5 & 5.24  & \textbf{58409.5} & \textbf{5.16 } & 59218.5 & 6.42  \\
			o2p2b2e2 & 57944 & 10.96  & \textbf{57551} & 3.44  & 57551 & 3.55  & 57634.5 & \textbf{2.57 } & 57899.5 & 3.97  & 57776 & 3.92  \\
			o3p1b1e1 & 55610 & 13.03  & 55262 & 4.09  & 55319 & 4.43  & 55605 & 4.59  & 55298 & 4.07  & \textbf{55153} & \textbf{3.99 } \\
			o3p1b2e1 & 54199 & 5.58  & 54025 & 2.48  & 54017 & \textbf{1.29 } & 53983 & 2.37  & \textbf{53973} & 2.53  & 54025 & 2.54  \\
			o3p1b1e2 & 55573 & 16.56  & 55399 & 5.18  & \textbf{55044.5} & \textbf{4.27 } & 55580.5 & 5.19  & 55739 & 5.86  & 55764.5 & 5.74  \\
			o3p1b2e2 & 54107.5 & 8.73  & 53869 & 2.69  & 53744 & 2.75  & \textbf{53732} & \textbf{2.54 } & 54102 & 2.88  & 53760 & 2.74  \\
			o2p2b1e2-t & 60808 & 32.62  & 59056 & 6.17  & \textbf{58566} & \textbf{5.41 } & 60011 & 7.66  & 58770.5 & 5.99  & 58822 & 5.96  \\
			o2p2b2e2-t & 57831 & 9.76  & 58078 & 3.84  & 57635 & \textbf{3.45 } & 57946 & 3.90  & \textbf{57537} & 3.49  & 57780 & 3.64  \\
			o3p1b1e1-t & 55408 & 11.78  & 55185 & 4.21  & 55215 & \textbf{3.94 } & 55249 & 4.00  & \textbf{55175} & 3.95  & 55423 & 4.42  \\
			o3p1b2e1-t & 54297 & 5.93  & 53970 & 1.67  & 53960 & \textbf{0.37 } & \textbf{53952} & 1.65  & 54109 & 1.76  & 54035 & 2.30  \\
			o3p1b1e2-t & 55868.5 & 19.78  & \textbf{55183} & \textbf{2.10 } & 55435.5 & 4.04  & 55502 & 4.98  & 55702.5 & 4.71  & 55301 & 4.37  \\
			o3p1b2e2-t & 54110 & 8.56  & 54363 & 3.62  & \textbf{53937} & \textbf{1.48 } & 54193 & 3.36  & 54167 & 2.28  & 54109 & 2.28  \\
			\bottomrule
			\end{tabular}%
		\begin{tablenotes}
\scriptsize
\item[1] Standard: Standard Benders decomposition (Sect. \ref{sec_bd_reform}); LP: Two-phase Benders decomposition (Sect. \ref{sec_imple}); LP-Pareto: Two-phase Benders decomposition with only Pareto-optimal cuts (Sect. \ref{sec_imple}); LP-Pareto-2: Two-phase Benders decomposition with both original and Pareto-optimal cuts (Sect. \ref{sec_imple}); LP-OC: Two-phase Benders decomposition with only Open-Close cuts (Sect. \ref{sec_imple}); LP-OC-2: Two-phase Benders decomposition with both original and Open-Close cuts (Sect. \ref{sec_imple}).
\end{tablenotes}
		\label{tab: small UB LB}%
	\end{table}%

	Firstly, it shows that in easy instances where the optimality can be reached, the standard setting of Benders decomposition requires the least computation times. This is due to the fact that the integer BMPs are easy to solve, and optimality cuts that are effective for the integer version of the problem can be directly and quickly generated. Secondly, when the two-phase approach (LP) is adopted, the lower bounds significantly improve in more difficult instances, which results in smaller optimality gaps. Thirdly, comparing the performances of the two strengthened optimality cuts, LP-Pareto achieved best optimality gaps in 9 out of 22 instances and LP-Pareto-2 achieved in 5. Further, the average optimality gaps of LP-Pareto and LP-Pareto-2 are 2.3\% and 2.59\%, respectively. While that of the LP-OC and LP-OC-2 are 2.71\% and 2.70\%, respectively. This shows the advantage of the Pareto-optimal cuts in both one- and two-cut versions. It is worth noting that the average optimality gap for LP is 2.52\%, which is better than the two Open-Close cut versions. This is due to the fact that the additional times to solve cut strengthening problems lead to less iterations performed, which results in worse performances. However, when using Pareto-optimal cuts, such disadvantage is remedied by the good quality of optimality cuts and better results are obtained. Still, LP-OC and LP-OC-2 produced best upper bounds in 6 instances, showcasing the benefits of Open-Close cuts. 

	Finally, considering the good performance of LP-Pareto in terms of both its smallest average optimality gap and its ability to find good upper bounds in most instances, we choose it as the basic algorithm setting for the following experiments.
	\subsubsection{Performance of the Benders decomposition algorithm.}\label{sec_performance}
	To see the correctness and performance of the proposed Benders decomposition algorithm, we also solve the instances using Gurobi and a callback approach. Due to the vast number of passenger paths which cannot be enumerated a priori, the Gurobi approach solves an alternative arc-based model shown in Appendix~\ref{app: arc formulation}, where passenger variables are arc-based. It is worth noting that the number of transfers could not be restricted in the arc-based model, and therefore this approach only solves instances without such requirement. Regarding the callback approach, we adopt the callback function of Gurobi to realize a single-tree implementation of Benders decomposition. The BMP is solved by Gurobi, and original and Pareto-optimal Benders optimality cuts are identified and generated for fractional and integer nodes of the search tree. A two-phase approach is also adopted here, which first solves the LP relaxation using standard Benders decomposition algorithm and then incorporates the generated cuts in the root node of the search tree. 

	Table \ref{tab: small benchmark} reports the computational results on the previously generated instances. It firstly shows that the arc-based model can be solved to optimality by Gurobi in very short times. This indicates that for small instances without transfer restrictions, solving the arc-based model by Gurobi is the best option. Second, our Benders decomposition algorithm (BD) outperforms the callback approach in terms of both smaller optimality gaps and less memory usage. The callback approach requires up to 70 GB of memory usage while our BD only requires memory usage of less than 12.3 GB. Finally, compared with the optimal solutions produced by Gurobi, the true qualities of the solutions obtained by BD can be analyzed. As shown by Gap$^2$(\%), the optimality gaps of BD are within 1.3\% among all the instances, which showcases the good performance of BD with respect to the quality of feasible solutions.  

	\begin{table}[htb] \tiny
		\centering
		\caption{Comparison results of different solution approaches}
		\begin{tabular}{cccccccccc}
			\toprule
			\multirow{2}[0]{*}{Instance} & \multicolumn{2}{c}{Gurobi} & \multicolumn{3}{c}{Callback} & \multicolumn{4}{c}{BD} \\
			\cmidrule(lr){2-3} \cmidrule(lr){4-6}\cmidrule(lr){7-10}
				  & UB    & CPU(s)/Gap(\%) & UB    & CPU(s)/Gap(\%) & Memory(GB) & UB    & CPU(s)/Gap(\%) & Memory(GB) & Gap$^2$(\%) \\
			\midrule
			o1p1b1e1 & 62007 & 1s    & 62007 & 15s   & 1.7   & 62007 & 12s   & 0.5   & 0.0  \\
			o1p1b2e1 & 61861 & 1s    & 61861 & 9s    & 0.6   & 61861 & 8s    & 0.3   & 0.0  \\
			o1p2b1e1 & 61096 & 1s    & 61096 & 336s  & 34    & 61149 & 2.1   & 6.2   & 0.1  \\
			o1p2b2e1 & 60666 & 1s    & 60666 & 148s  & 3     & 60666 & 76s   & 2.2   & 0.0  \\
			o2p1b1e2 & 59579 & 2s    & 60546.4 & 2.6   & 47    & 59645 & 1.0   & 6.1   & 0.1  \\
			o2p1b2e2 & 59027 & 1s    & 59027 & 159s  & 17    & 59027 & 168s  & 8.7   & 0.0  \\
			o2p1b1e1 & 59727 & 1s    & 59727 & 592s  & 6     & 59727 & 0.1   & 6.3   & 0.0  \\
			o2p1b2e1 & 59364 & 1s    & 59364 & 93s   & 10    & 59364 & 184s  & 5     & 0.0  \\
			o2p2b1e1 & 58558 & 2s    & 67043.9 & 18.0  & 60    & 58912 & 4.9   & 5.9   & 0.6  \\
			o2p2b2e1 & 58215.5 & 1s    & 61093.9 & 9.4   & 76    & 58275.5 & 1.8   & 5.5   & 0.1  \\
			o2p2b1e2 & 58050 & 4s    & 63165.9 & 14.9  & 65    & 58791 & 5.7   & 10.1  & 1.3  \\
			o2p2b2e2 & 57339 & 1s    & 61356.4 & 12.1  & 70    & 57551 & 3.5   & 8.1   & 0.4  \\
			o3p1b1e1 & 54981 & 4s    & 60642.9 & 14.4  & 64    & 55319 & 4.4   & 5.4   & 0.6  \\
			o3p1b2e1 & 53829 & 1s    & 61861 & 15.9  & 54    & 54017 & 1.3   & 5.1   & 0.3  \\
			o3p1b1e2 & 54863 & 7s    & 61411.4 & 16.3  & 48    & 55044.5 & 4.3   & 6.7   & 0.3  \\
			o3p1b2e2 & 53570 & 2s    & 57220.9 & 10.5  & 64    & 53744 & 2.7   & 11.9  & 0.3  \\
			o2p2b1e2-t & -     & -     & 67422.9 & 20.1  & 62    & 58566 & 5.4   & 10.7  & - \\
			o2p2b2e2-t & -     & -     & 59548.9 & 9.1   & 70    & 57635 & 3.4   & 8.4   & - \\
			o3p1b1e1-t & -     & -     & 61411.9 & 15.0  & 60    & 55215 & 3.9   & 7.6   & - \\
			o3p1b2e1-t & -     & -     & 54394.9 & 3.9   & 60    & 53960 & 0.4   & 12.3  & - \\
			o3p1b1e2-t & -     & -     & 57883.9 & 10.7  & 51.3  & 55435.5 & 4.0   & 6.6   & - \\
			o3p1b2e2-t & -     & -     & 56999.4 & 10.0  & 56.6  & 53937 & 1.5   & 11.6  & - \\
			\bottomrule
			\end{tabular}%
			\begin{tablenotes}
\scriptsize
\item[1] Callback: The callback approach using Gurobi where Pareto-optimal Benders optimality cuts are added for integer and fractional nodes in the search tree (Sect.\ref{sec_performance}; BD: The proposed Benders decomposition algorithm using two-phase decomposition and only Pareto-optimal Benders optimality cuts (Sect. \ref{sec:algo}).
\end{tablenotes}
		\label{tab: small benchmark}%
	\end{table}%

	\subsubsection{Effects of stop-skipping technique of extra train services.}
\begin{figure}[htb]
		\centering
		\subfigure[Objective]{
		\includegraphics[width=0.31\textwidth]{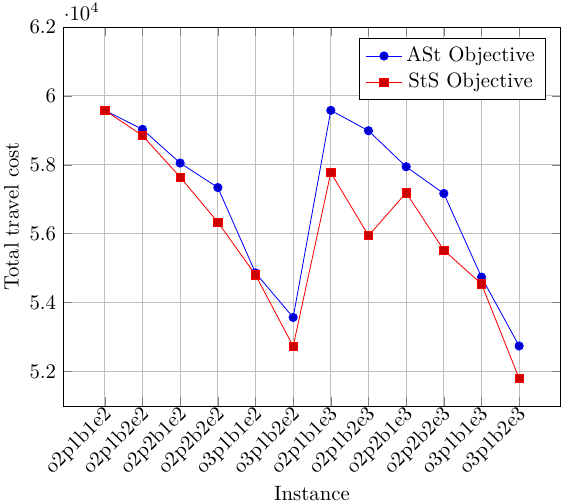}\label{fig: skip-stop-obj}
		}
		\subfigure[Average in-vehicle cost]{
		\includegraphics[width=0.31\textwidth]{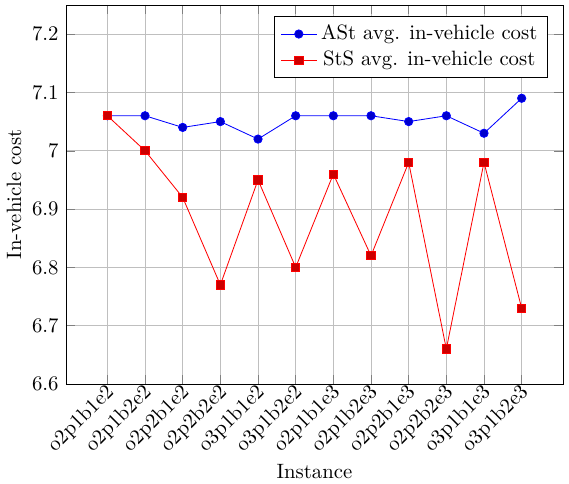}\label{fig: skip-stop-iv}
		}
		\subfigure[Computation time]{
		\includegraphics[width=0.31\textwidth]{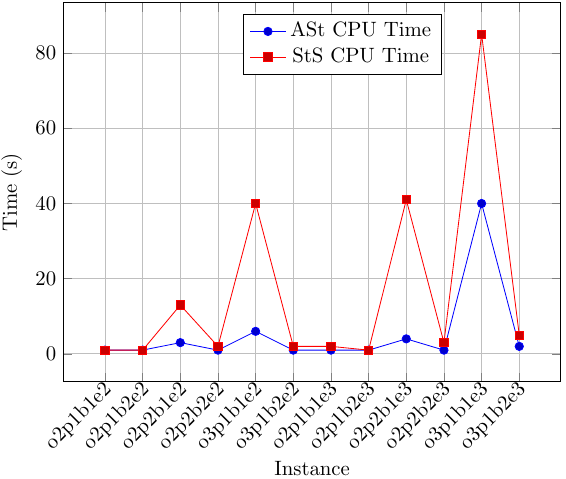}\label{fig: skip-stop-cpu}
		}
	 	\caption{\centering{Experimental results of the stop-skipping technique of extra train services}}\label{fig: skip-stop}
		{}
	\end{figure}

	In this experiment, we test the effect of the stop-skipping technique of extra train services. As introduced before, the stop plans of extra train services can be implied by the selection of arcs in the time-space network. Figure~\ref{fig: skip-stop} shows the comparison results of all-stop (ASt) and stop-skipping (StS) techniques in 12 instances. In the former, only train arcs associated with station departure/arrival stop nodes are included in the subnetwork of extra trains. For accurate analysis, all instances are solved to optimality by Gurobi. 

	Figure \ref{fig: skip-stop-obj} shows that utilizing the stop-skipping technique can reduce passenger travel times across all instances, which is mainly due to the decrease in the in-vehicle times of passengers indicated by Figure \ref{fig: skip-stop-iv}. This is because some stations are skipped by extra trains, which saves the station dwell times and additional section running times for the passing passengers. However, instances allowing the stop-skipping technique are more difficult to solve and require more computation times as indicated by Figure \ref{fig: skip-stop-cpu}. In instance o2p2b1e3, up to ten times of ASt computation times are needed by StS. Therefore, the test of stop-skipping technique is limited to the Toy case. 
	
	\subsection{Large case study}
	In this section, we test the algorithm on the Lower Saxony case in LinTim, which contains 34 stations and 70 sections. 836 time-dependent passengers are considered in this case. A one-hour periodic timetable with 20 train services is given. Two types of rolling stock units with 200 and 300 seats are considered. The time discretization is 2 minutes and we consider the scheduling time horizon of 6 hours. Other details of the case can be seen in Appendix~\ref{app: lower saxony}.
	
	The periodicity level of original trains $\xi$ is set as 0.8. The LP phase terminates when the optimality gap is within 5\%. $\alpha_0, \bar \alpha$\oxc and $\epsilon$ are set as 7\%, 3\%\oxc and 3\%. The computation time limit is set as 6 hours. We also consider different problem settings for the case. Three levels of original timetable deviations are considered, where o1, o2\oxc and o3 represent the deviation time windows of 0, 4\oxc and 8 minutes, respectively. Two time windows for extra trains are considered, where e1 and e2 represent the deviation of 4 and 8 minutes, respectively. Besides, e0 represents no extra trains are allowed to be inserted. Further, two levels of transfer restrictions are considered, where t1 and t2 represents that 1 and 2 transfers are allowed in each passenger route.
	\subsubsection{Network preprocessing.}
	To see the effects of our network preprocessing techniques, we record the numbers of original train and passenger arc variables if different preprocessing techniques are applied, and the results are reported in Table \ref{tab: prep}. It shows that the number of original train arcs (\# Train Arc) slightly reduces if we truncate the network considering the right time boundary of the planning horizon. The number of passenger arcs (\# Pax Arc) reduces by 16.9\% if the network is truncated similarly. If only the arcs that might be used by trains are kept in passengers' subnetworks, the number of passenger arcs reduces by roughly 30\% under the three settings of deviation levels. This demonstrates that our network preprocessing techniques can significantly reduce the size of the problem. While it is noted that such reduction becomes less effective when the deviation times enlarge, as more arcs are possibly used by trains, leading to less deleted passenger arcs. 
	\begin{table}[htb] \tiny
		\centering
		\caption{Experiment results of network preprocessing techniques}
		\begin{tabular}{ccccccc}
			\toprule
			\multirow{2}[0]{*}{Deviation level} & \multicolumn{2}{c}{No Prep.} & Train Prep. & Pax Prep. & Train Prep. + Pax Prep. + Active Arc & Pax Arc Reduction \\
			\cmidrule(lr){2-3} \cmidrule(lr){4-4}\cmidrule(lr){5-5}\cmidrule(lr){6-6}\cmidrule(lr){7-7}
				  & \# Train Arc & \# Pax Arc & \# Train Arc & \# Pax Arc & \# Pax Arc & (\%) \\
			\midrule
			o1    & 2216  & 2096099 & 2216  & 1740800 & 1392024 & 33.6  \\
			o2    & 6640  & 2096099 & 6581  & 1740800 & 1459058 & 30.4  \\
			o3    & 11058 & 2096099 & 10935 & 1740800 & 1510935 & 27.9  \\
			\bottomrule
			\end{tabular}%
\begin{tablenotes}
\scriptsize
\item[1] No Prep.: No network preprocessing performed; Train Prep.: Train network preprocessing by operation ending time (Sect. \ref{sec_prep}); Pax Prep.: Passenger network preprocessing by latest arrival time (Sect. \ref{sec_prep}); Active arc: Passenger network preprocessing considering active train arcs (Sect. \ref{sec_prep}).
\end{tablenotes}
		\label{tab: prep}%
	\end{table}%

	\subsubsection{Acceleration of column generation.}
	To see the performance of acceleration techniques for column generation, we conduct experiments with different time windows for original and extra trains, and also the maximum number of transfers in passengers' routes, which all influence the column structures. Table \ref{tab: colgen acc} reports the average computation times (Avg. Time(s)) of initial ten iterations of column generation and the total number of columns generated (\# Col.). ``Standard'' represents the implementation without acceleration techniques, ``PP'' and ``CC'' indicate the implementations of partial pricing and column copying techniques, respectively. ``CC + PP'' represents the hybrid implementation of the two. Firstly, it shows that as the time windows of existing and extra trains enlarge, the number of columns generated increases. Besides, the number of columns also increases significantly with the maximum number of transfers in a passenger's route, since more paths are considered feasible. Second, PP decreases the number of generated columns, which indicates that many repetitive columns can be avoided, while CC might increase the number of generated columns, as many unfavorable columns might be copied in the first place. Both implementations obtained shortest computation times in two instances. Third, the hybrid implementation obtained the shortest computation times in half of the instances, and the PP technique managed to reduce the number of columns which might be enlarged by CC initially, indicating that the hybrid implementation can accelerate computation and also reduce the number of generated columns. In difficult instances o2e1t2 and o2e2t2, the time and column savings can reach more than ten and eight percent, respectively.

	\begin{table}[htb] \scriptsize
		\centering
		\caption{Experiment results of acceleration techniques for column generation}
		\begin{tabular}{ccccccccc}
			\toprule
			\multirow{2}[0]{*}{Instance} & \multicolumn{2}{c}{Standard} & \multicolumn{2}{c}{PP} & \multicolumn{2}{c}{CC} & \multicolumn{2}{c}{CC+PP} \\
			\cmidrule(lr){2-3} \cmidrule(lr){4-5}\cmidrule(lr){6-7}\cmidrule(lr){8-9}
				  & Avg. Time(s) & \# Col. & Avg. Time(s)  & \# Col. & Avg. Time(s) & \# Col. & Avg. Time(s)  & \# Col. \\
			\midrule
			o1e1t1 & 34.4  & 25505 & 32.8  & 23770 & 34.8  & 26062 & \textbf{30.8} & 23807 \\
			o1e2t1 & 39.6  & 29098 & 36.4  & 27271 & \textbf{33.1} & 27907 & 33.8  & 27055 \\
			o2e1t1 & 46.6  & 39359 & \textbf{43.3} & 35462 & 47.5  & 39850 & 45.6  & 37350 \\
			o2e2t1 & 49.9  & 39689 & 50.9  & 38983 & 50.1  & 42761 & \textbf{48.2} & 40333 \\
			o1e1t2 & 48.5  & 35991 & 47.5  & 34302 & \textbf{42.8} & 37014 & 45.5  & 33119 \\
			o1e2t2 & 60.3  & 43639 & \textbf{58.6} & 41124 & 59.9  & 43720 & 63.7  & 41381 \\
			o2e1t2 & 88.5  & 59377 & 79.8  & 52032 & 80.8  & 61372 & \textbf{78.9} & 54135 \\
			o2e2t2 & 83.5  & 60396 & 83.8  & 54995 & 79.9  & 62373 & \textbf{72.6} & 54196 \\
			\bottomrule
			\end{tabular}%
\begin{tablenotes}
\scriptsize
\item[1] Standard: Standard column generation (Sect. \ref{sec_cg}); PP: Partial pricing technique (Sect. \ref{sec_imple}); CC: Column-copying technique (Sect. \ref{sec_imple}); CC+PP: Combined partial pricing and column-copying technique; Bold numbers indicate the shortest computation times of each instance.
\end{tablenotes}
		\label{tab: colgen acc}%
	\end{table}%

	\subsubsection{Benefits of hybrid periodicity.}

\begin{figure}[htb]
		\centering
		\subfigure[Upper, lower bounds and unrouted costs]{
		\includegraphics[width=0.4\textwidth]{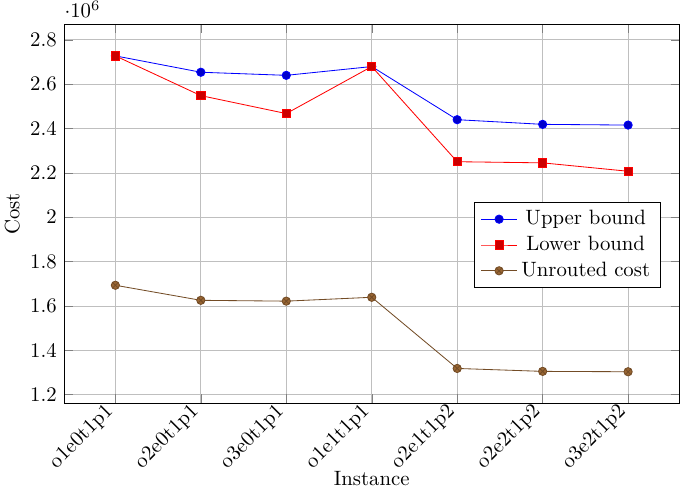}\label{fig_hybrid_a}
		}
		\quad
		\subfigure[Travel and transfer costs]{
		\includegraphics[width=0.4\textwidth]{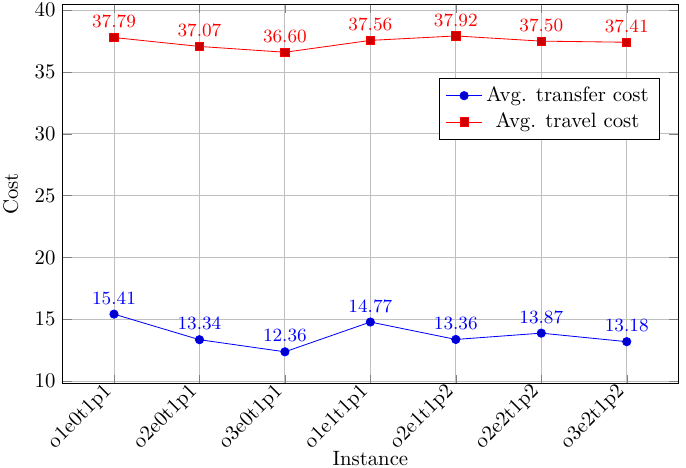}\label{fig_hybrid_b}
		}
	 	\caption{\centering{Experiment results with different periodicity levels and time windows}}\label{fig: hybrid periodicity}
	\end{figure}

	To see the benefits of the hybrid periodicity of timetables, experiments with different periodicity levels and allowable time deviations are conducted for the Lower Saxony case. All instances are solved by the proposed two-phase Benders decomposition algorithm with only Pareto-optimal optimality cuts as introduced in Section \ref{sec_imple}. Figure \ref{fig_hybrid_a} shows the upper, lower bounds and passenger unrouted costs, and Figure \ref{fig_hybrid_b} shows the average transfer and travel costs of each passenger. 
	
	Firstly, it shows that when all original trains are kept and their timetables are fixed (instances with o1p1), optimal solutions can be obtained. While the problem becomes much more difficult for other instances when the flexibility is allowed for original trains. For the most difficult instance o3e2t1p2, the optimality gap reaches 8.6\%. Second, when all the original trains are kept and no extra trains are allowed to be inserted (instances with e0p1), the qualities of solutions slightly improve with the allowable time windows of original trains, which corresponds with the decrease of travel and unrouted costs of passengers. Third, the solution quality of o1e1t1p1 is better than that of o1e0t1p1, which is because two extra trains are scheduled in the former one (while no original trains are canceled). The utilized operation budgets in the two instances both equal the given maximum budget, which indicates that the extra trains help circulate the RSUs in a budget-efficient way, which allows two extra trains to be inserted by changing the types of RSUs assigned to each original train. Finally, it is shown that when original trains can be canceled (instances with p2), a sharp decrease in passenger travel costs is observed, which corresponds with the sharp decrease in passenger unrouted cost. This is due to the fact that the cancellation of original train services in non-peak periods saves the operation budget for operating the extra train services in the peak ones, thus reducing the number of unrouted passengers and the unrouted costs significantly. 

	Therefore, the experiment shows that the hybrid periodicity can reduce travel costs for routed passengers by modifying the existing train schedules, and also help better allocate limited transport capacity to better satisfy peak-hour passenger demands. Comparing with the original timetable o1e0t1p1, the most flexible instance o3e2t1p2 reduces the total passenger cost by 11.4\%. 

	\subsubsection{Routing a subset of passenger groups.}
	To test the effect of routing a subset of passenger groups using HPTTP-PSR, we conduct experiments on three larger instances with the time horizon of 12 hours. Such a longer  time horizon thus involves more passenger groups, which further increases the computation burden of the Benders subproblem. The three instances share the same problem setting except the passenger demand, where instance ``Standard'' has the basic passenger demand similar to the previous section, while in ``Random'' and ``High'', the demand is randomly perturbed and increased, respectively. We consider four solution approaches, where HPTTP is the two-phase Benders decomposition algorithm with only Pareto optimal cuts for the original model introduced in Section \ref{sec_imple}, HPTTP-PSR is the Benders decomposition algorithm for the passenger subset routing model introduced in Section \ref{sec: subset routing}, HPTTP-Delete is to first delete several passenger groups before solving HPTTP but reroute all passenger groups on the resulting timetable\oxc and None-Routed is a special case when all passenger groups are deleted in HPTTP-Delete. By experiments, we choose to fix or delete 400 passenger groups with the least number of passengers.

	Table \ref{tab: psr} shows the computation results. First, it shows that HPTTP-PSR obtained best objectives among all the instances, indicating the effectiveness of the approach. Second, HPTTP-Delete performs better than HPTTP in the first instance, but worse in the rest. This shows that when completely disregarding a subset of passenger groups within optimization, the results can be worse. Also, when none of the passengers are routed, the results are much worse. Together, it indicates that by taking into account a subset of passenger groups partially, i.e., by assigning fixed routes for them, better overall results can be obtained. Second, column ``BMP time proportion (\%)'' indicates that BMP-PSR is much more difficult to solve than BMP due to the new variables and constraints w.r.t. fixed passengers. Finally, HPTTP-PSR generates less passenger columns than the other approaches, which is due to the fact that fixed passengers are not considered in the subproblem, and also less computation times are devoted to the subproblem computation. 

	\begin{table}[htb] 
		\scriptsize
		\centering
		\caption{Experiment results of passenger subset routing test}
		\begin{tabular}{ccccc}
			\toprule
			Case  & Method & Objective & BMP time proportion (\%) & \#Pax Col \\
			\midrule
    \multirow{4}[0]{*}{Standard} & HPTTP & 4547380 & 5.4   & 291561 \\
          & HPTTP-Delete & 4495360 & 16.9  & 306658 \\
          & HPTTP-PSR & 4492890 & 28.4  & 282440 \\
          & None routed & 4861360 & -     & - \\
		  \midrule
    \multirow{4}[0]{*}{Random} & HPTTP & 4531640 & 4.5   & 296178 \\
          & HPTTP-Delete & 4578310 & 7.7   & 281076 \\
          & HPTTP-PSR & 4527500 & 21    & 275538 \\
          & None routed & 4825710 & -     & - \\
		  \midrule
    \multirow{4}[0]{*}{High} & HPTTP & 5630180 & 3.2   & 273210 \\
          & HPTTP-Delete & 5637960 & 4.5   & 264737 \\
          & HPTTP-PSR & 5592580 & 10.7  & 257211 \\
          & None routed & 5964630 & -     & - \\
			\bottomrule
			\end{tabular}%
			\begin{tablenotes}
\scriptsize
\item[1] HPTTP: Two-phase Benders decomposition algorithm with only Pareto-optimal cuts for the initial model HPTTP (Sect. \ref{sec_imple}); HPTTP-PSR: Benders decomposition algorithm for passenger subset routing model HPTTP-PSR (Sect. \ref{sec: subset routing}); HPTTP-Delete: First delete several passenger groups before solving HPTTP, then reroute all passengers on the resulting timetable; None-routed: A special case of HPTTP-Delete where all passengers are deleted.
\end{tablenotes}
		\label{tab: psr}%
	\end{table}%

\section{Conclusions} \label{sec:conclu}
This study considered a hybrid periodic train timetabling problem which aims to increase the flexibility of a given periodic timetable by schedule adjustment, extra train scheduling\oxc and efficient rolling stock circulation. Based on a time-space network representation of the problem, a dynamic service network design model was developed, which was then solved by a decomposition algorithm integrating Benders decomposition and column generation. Various techniques were proposed to accelerate both the Benders master and subproblem. Further, we investigated a problem variant which fixed part of passenger groups and only routed a subset of passenger groups. 
Numerical experiments on artificial small instances and also the Lower Saxony case in LinTim showcased the good performance of the algorithm and also the effectiveness of various acceleration techniques. Compared with the original periodic timetable, the timetable with hybrid periodicity managed to reduce the passenger travel costs by 11.4\%, showcasing the benefits of the hybrid periodicity. Besides, the problem variant that routes a subset of passenger groups also achieved better results in large problem instances.

There are several possible extensions to this study. First, the coupling and decoupling of RSUs can be considered to further increase the flexibility of rolling stock circulation and better match passenger demand. Second, it remains to consider different passenger demands in different days, and develop corresponding hybrid periodic timetables and also consistent rolling stock circulations for a longer time horizon.


\begingroup
\tiny
\let\oldbibitem\bibitem
\renewcommand{\bibitem}{\setlength{\itemsep}{0pt}\oldbibitem}
\setlength{\bibsep}{0pt}
\bibliographystyle{informs2014trsc} 
\bibliography{cas-refs.bib} 
\endgroup

\begin{APPENDICES}
\newpage
\section{Tables of notations}\label{app: tables}
Table \ref{tab:model sets and indexes} summarizes the notations related to the problem description and Table \ref{tab:ts net sets and indexes} summarizes the notations related to the time-space network representation of the problem.
\begin{table}[htbp] \footnotesize
	\caption{Notations used in the problem description}
	\begin{tabular*}{\textwidth}{p{0.11\textwidth}p{0.08\textwidth}p{0.81\textwidth}}\toprule
            \textbf{Component} & \textbf{Type} & \textbf{Definition} \\
		\midrule
		$\mathcal{M}$ & Set & Stations, indexed by $m$\\
		$\mathcal{M}^\textup{ter}$ & Set & Terminal stations where trains can originate and terminate \\
            $\mathcal{M}_m$ & Set& Adjacent stations of station $m$\\
		$\mathcal{E}$ & Set & Railway sections, indexed by $e$ \\
		$\mathcal{T}$ & Set & Discrete timestamps in the time horizon of train timetable, indexed by $t$\\
		$\mathcal{T}^\textup{inv}$ & Set & Discrete timestamps in the time horizon of rolling stock circulation\\
		$\mathcal{T}_{r}^\textup{all}$ & Set & Discrete timestamps in the available departure time window of passenger group $r$\\
		${\cal T}_r^\textup{pre}$ & Set & Discrete timestamps in the preferred departure time window of passenger group $r$\\
		$\mathcal{K}$ & Set & Original trains, indexed by $k$\\
		$\mathcal{L}$ & Set & Original lines, indexed by $l$\\
		$\mathcal{K}_l$ & Set & Trains belonging to line $l$\\
		$\mathcal{R}$ & Set & Passenger groups, indexed by $r$ \\
		$\mathcal{U}$ & Set & RSU types, indexed by $u$ \\
            \midrule
		$m_r^\textup{org}, m_r^\textup{des}$ & Index & Origin and destination station of passenger group $r$\\    
		\midrule
		$\xi$ & Param. & Minimum proportion of the operated trains in each periodic line\\
		$\tau$ & Param.& Allowed time deviation of original trains' departure and arrival times\\
		$t_r^\textup{min}, t_r^\textup{max}$ & Param.& Smallest and biggest timestamps in the preferred departure times of passenger group $r$\\
		$\rho^\textup{conn}$ & Param.& Minimum transition time of RSUs\\
		$\varpi _m^u$ & Param.& Initial number of RSUs of type $u$ at station $m$\\
		$B$ & Param.& Budget for seat-kilometers of trains being operated\\
		$f_r$ & Param.& Penalty cost for each unsatisfied passenger in passenger group $r$ \\
		$p_u$ & Param.& Seat capacity of RSU of type $u$ \\
		$g_r$ & Param.& Number of passengers in passenger group $r$ \\  
		$\varUpsilon$ & Param.& Maximum number of transfers in a passenger route \\
		\bottomrule
	\end{tabular*}
	\label{tab:model sets and indexes}
\end{table}

\newpage

\begin{table}[htbp] \footnotesize
	\centering
	\caption{Notations used in the time-space network}
	\begin{tabular*}{\textwidth}{p{0.15\textwidth}p{0.07\textwidth}p{0.78\textwidth}}\toprule
		\textbf{Component} & \textbf{Type} & \textbf{Definition} \\
		\midrule
		$\cal{N}$ & Set & Physical nodes, indexed by $n$ \\
		${\cal{N}}_k$ & Set & Physical nodes of train $k$, ${\cal N}_k \subseteq {\cal N}$ \\
		${\cal{N}}^\textup{inv}$ & Set & Station inventory nodes \\
		${\cal{N}}^\textup{od}$ & Set & Passenger origin and destination nodes \\
		$\mathcal{V}$  & Set & Time-space vertices, indexed by $v$ \\
		$\mathcal{V}_k$ & Set & Time-space vertices of train $k$ \\
		$\mathcal{V}^\textup{inv}$ & Set & Time-space inventory vertices \\
		$\mathcal{V}^\textup{ex}$ & Set & Time-space vertices that can be visited by extra trains \\
		$\mathcal{V}^\textup{bound}$ & Set & Time-space vertices at the station boundary nodes \\
		$\mathcal{A}$ & Set & Time-space arcs, indexed by $a$ \\
		${\cal{A}}_k$ & Set & Time-space arcs of train $k$ \\
		${\cal{A}}^\textup{ex}$ & Set & Time-space arcs that can be used by extra trains\\
		${\cal{A}}^\textup{vir}_k,{\cal{A}}^{\textup{ex},\textup{vir}}$ & Set & Virtual arcs of train $k$ and extra trains\\
		${\cal{A}}^\textup{tr}$ & Set & Non-virtual train time-space arcs \\
		${\cal{A}}^\textup{sec}$ & Set & Train section arcs \\
    		${\cal{A}}^\textup{inv}$ & Set & Station inventory arcs of RSU of type $u$ at stations, ${\cal A}^\textup{inv} \subseteq {\cal A}$ \\
		$\phi_{n,t}$ & Set & Incompatible arcs associated with time-space vertex $(n,t)$\\
            \midrule
		$o_k, d_k$ & Index & Starting and ending station inventory nodes of train $k$\\
		$o_r, d_r$ & Index & Origin and destination nodes of passenger group $r$\\
		$d_m^{m',\textup{stop}}, a_m^{m', \textup{stop}}$ & Index & Departure and arrival stop node of station $m$ associated with the adjacent station $m'$\\
		$d_m^{m',\textup{skip}}, a_m^{m', \textup{skip}}$ & Index & Departure and arrival skip node of station $m$ associated with the adjacent station $m'$\\
		$n^\textup{trans}_m$ & Index & Passenger transfer node at station $m$ \\
		$n^\textup{inv}_m$ & Index & Rolling stock unit inventory node at station $m$ \\
            \midrule
		$w^\textup{shift}, w^\textup{wait}, w^\textup{veh}$, $w^\textup{trans}$ &Param. & Cost coefficients of passenger departure shift, waiting, in-vehicle and transfer times. \\
		$\rho^\textup{sec}_e$ &Param. & Minimum running time of section $e$ \\
		$\rho^\textup{acc}_e, \rho^\textup{dec}_e$ &Param. & Additional running times due to train acceleration and deceleration at section $e$\\
		$\underline{\rho}^\textup{dwell}_m, \overline{\rho}^\textup{dwell}_m$ &Param. & Minimum and maximum train dwell times at stations\\
		$\rho^\textup{trans}_m$ &Param. & Minimum passenger transfer time at station $m$\\
		$\beta_a$ &Param. & Distance in kilometers of section arc $a$\\
		$\beta_k$ &Param. & Distance in kilometers of original train $k$\\
		$c_a^r$ &Param. & Cost of arc $a$ for passenger group $r$\\
		\bottomrule
	\end{tabular*}
	\label{tab:ts net sets and indexes}
\end{table}

\section{Transfer capturing capability of the time-space network}
\label{app: transfer capture}
	It is noted that by introducing transfer arcs in the time-space network, passenger transfers can be considered in such representations \citep{Zhan2021-ADMM,Xu2021-LR,Yao2023-ADMM}. However, transfer costs might not be captured if no additional transfer nodes are introduced, which is the case in existing studies.

	Usually, passenger transfer arcs are set to link the arrival and departure vertices of a station to represent the passenger walking process in transfers, and additional waiting activity can be represented by using the unit waiting arcs at the departure nodes. However, in this case, the train dwell arc and passenger transfer arc link the same physical nodes, which might result in the miscalculation of passenger transfer costs. Figure \ref{fig: transfer capture} gives a simple example, where passengers are transferring from train 1 to train 2. If no transfer nodes are introduced, Figure \ref{fig: transfer capture}(a) shows the correct transfer route of passengers: they disembark from train 1, use the transfer arc to move to the departure node\oxc and wait for three time intervals before boarding train 2. However, as Figure \ref{fig: transfer capture}(b) illustrates, passengers can also use the dwell arc of train 1 to move to the departure node and wait for one time interval before boarding train 2. The second route is less costly as the transfer arcs are much more costly than the train dwell arcs, and therefore route 2 is selected in the results as it dominates route 1, leading to undetected transfer costs.
	\begin{figure}[htb]
		\centering
		\FIGURE
		{\includegraphics[width=\textwidth]{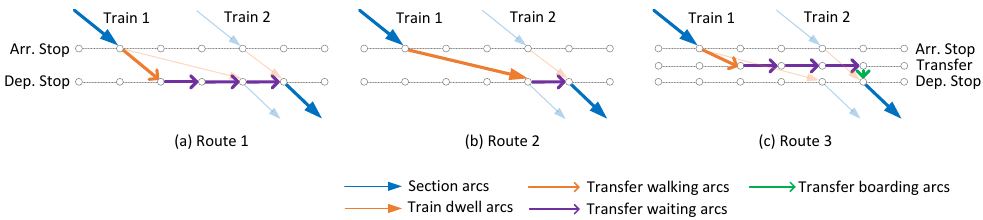}}
		{\centering{Illustration of transfer missing phenomenon}\label{fig: transfer capture}}
		{}
	\end{figure}

	To prevent this phenomenon, we introduce a transfer node at each station as shown in Figure \ref{fig: transfer capture}(c). The waiting arcs now link two adjacent vertices at the transfer node instead of the departure node, and passengers can only wait at the transfer node for their connecting train. To use the unit waiting arcs, passengers must use the transfer walking arc first to access the waiting arcs, and therefore transfer costs are captured.

\section{Different minimum headway times and incompatible arc sets}
\phantomsection
\label{app: headway}

Seven headway times are considered based on the stop and skip patterns of adjacent trains, which are shown in Figure \ref{fig: headway types}. For example, when the former train departs from the station after a dwelling operation and the second train passes the station, then the departure headway time $h_\textup{dp}$ should be imposed on the two trains' departure times, where ``d'' denotes ``departure'' and ``p'' denotes ``passing''. In the expressions of arrival headways, ``a'' denotes ``arrival''. Therefore, four departure headway times $h_\textup{dp}, h_\textup{dd}, h_\textup{pd}, h_\textup{pp}$ and three arrival headway times $h_\textup{ap},h_\textup{pa},h_\textup{aa}$ are considered. 

\begin{figure}[H]
	\centering
	\FIGURE
	{\includegraphics[width=\textwidth]{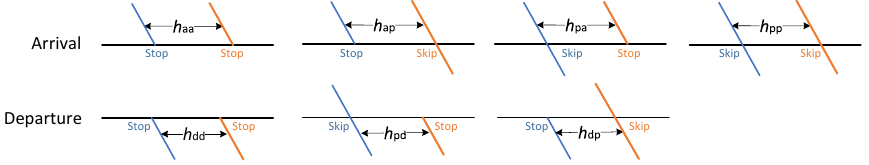}}
	{\centering{Illustration of different minimum headway times}\label{fig: headway types}}
	{}
\end{figure}

To model track capacity/headway requirements, we view each time-space vertex at departure and arrival nodes as a resource. Since nodes $d_m^{m', \textup{skip}}$ and $d_m^{m', \textup{stop}}$ both represent train departures bound for station $m'$ at station $m$, the section arcs then occupy the vertices at both nodes. This also applies to the two arrival nodes associated with the same adjacent station as well. Figure \ref{fig: vertex occup} gives an example of how section arcs occupy the vertices. Take the left section arc for example, it departs at time $t$ from skip node $d_{A}^{B, \textup{skip}}$, then it occupies vertices $\SetCond{(d_A^{B,\textup{skip}},t')}{t' \in {\cal T}:t \le t' \le t + {h_\textup{pp}} - 1}$ at the skip node, thus preventing the next skip departure in ${h_\textup{pp}}$ minutes. Besides, it also occupies vertices $\SetCond{(d_A^{B,\textup{stop}},t')}{t' \in {\cal T}: t \le t' \le t + {h_\textup{pd}} - 1}$ at the departure stop node, preventing the next stop departure in $h_\textup{pd}$ minutes.
\begin{figure}[htb]
	\centering
	\FIGURE
	{\includegraphics[width=0.7\textwidth]{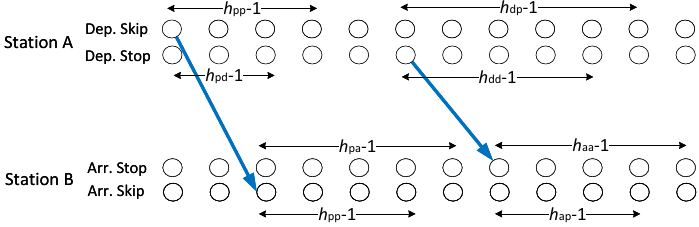}}
	{\centering{Illustration of vertex occupation of train section arcs}\label{fig: vertex occup}}
	{}
\end{figure}

When a vertex is occupied more than once, then the headway time ranges of the associated selected train arcs overlap at the timestamp of the vertex, meaning that the trains are not sufficiently spaced in time and the headway requirement is violated. Therefore, denote by $\phi_{n,t}$ the incompatible arc set of vertex $(n,t)$ containing all the arcs that occupy it. By requiring no more than one arc in $\phi_{n,t}$ can be selected, the headway requirement with respect to vertex $(n,t)$ is thus respected.

Eq. (\ref{eq: arc set}) shows the incompatible section arc sets of vertices in ${\cal V}^\textup{bound}$.

\vspace{-5pt}
\begingroup
\small
\begin{equation}
	\begin{split}
		{\phi _{d_m^{m',\textup{stop}},t}} &= {\rm{ }} \SetCond{(d_m^{m',\textup{stop}},n,t',t'') \in {\cal A}^\textup{sec}}{t' \in {\cal T} \cap [t - {h_\textup{dd}} + 1,t]}
		\cup \SetCond{(d_m^{m',\textup{skip}},n,t',t'') \in {\cal A}^\textup{sec}}{t' \in {\cal T} \cap [t - {h_\textup{pd}} + 1, t]}\\
		{\phi _{d_m^{m',\textup{skip}},t}} &= {\rm{ }} \SetCond{(d_m^{m',\textup{skip}},n,t',t'') \in {\cal A}^\textup{sec}}{t' \in {\cal T} \cap [t - {h_\textup{pp}} + 1 , t]} \cup \SetCond{(d_m^{m',\textup{stop}},n,t',t'') \in {\cal A}^\textup{sec}}{t' \in {\cal T} \cap [ t - {h_\textup{dp}} + 1,t]}\\
		{\phi _{a_m^{m',\textup{stop}},t}} &= \SetCond{(n,a_m^{m',\textup{stop}},t',t'') \in {\cal A}^\textup{sec}}{t'' \in {\cal T} \cap[ t - {h_\textup{aa}} + 1 ,t]} \cup \SetCond{(n,a_m^{m',\textup{skip}},t',t'') \in {\cal A}^\textup{sec}}{t'' \in {\cal T} \cap [t - {h_\textup{pa}} + 1 ,t]} \\
		{\phi _{a_m^{m',\textup{skip}},t}} &= \SetCond{(n,a_m^{m',\textup{skip}},t',t'') \in {\cal A}^\textup{sec}}{t'' \in {\cal T} \cap [t - {h_\textup{pp}} + 1 ,t]}  \cup \SetCond{(n,a_m^{m',\textup{stop}},t',t'') \in {\cal A}^\textup{sec}}{t'' \in {\cal T} \cap [t - {h_\textup{ap}} + 1,t]}
	\end{split}
	\label{eq: arc set}
\end{equation}
\endgroup

\section{Column generation for Open-Close cut strengthening model}\label{app: oc cut}

Here we detail how to use column generation to solve the cut strengthening model of the Open-Close cut and a simplification technique. Denote the dual variables of constraints (\ref{con:rc-DBSP-strengthen}) as $z_p^r$, we can write the dual of the problem:
\vspace{-5pt}
\begingroup
\small
\begin{align}
	\min \quad & \left( {c_p^r - \sum\limits_{a \in {{\cal O}^\textup{tr}}} {d_p^a\lambda _a^*}  - \mu _r^*} \right)z_p^r\\
	\textrm{s.t.} \quad &\sum\limits_{r \in {\cal R}} {\sum\limits_{p \in {{\cal P}_r}} {d_p^a} } z_p^r \le 1 \quad \forall a \in {\cal C}^\textup{tr}
	\label{con:str-packing}\\
	&z_p^r \ge 0 \quad \forall r \in {\cal R}, p \in {\cal P}_r
\end{align}
\endgroup
It is noted that this problem can also be solved by column generation, and the reduced cost of a path can be written as Eq. (\ref{rc:strengthen}). It can also be solved by the shortest path algorithm introduced earlier, but with the modification that the fixed optimal arc dual prices $\lambda_a^*$ are used on the open arcs, and the arc dual prices $\lambda_a$ of constraints (\ref{con:str-packing}) are used on the closed arcs. When the problem is solved, the additional Benders optimality cut can be generated.

\vspace{-10pt}
\begingroup
\small
\begin{equation}
	\begin{split}
	\textup{rc}_p^r = \left( {c_p^r - \sum\limits_{a \in {{\cal O}^\textup{tr}}} {d_p^a\lambda _a^*}  - \mu _r^*} \right) - &\sum\limits_{a \in {{\cal C}^\textup{tr}}} {d_p^a{\lambda _{a}}}  = c_p^r - \sum\limits_{a \in {{\cal O}^\textup{tr}}} {d_p^a\lambda _a^*}  - \sum\limits_{a \in {{\cal C}^\textup{tr}}} {d_p^a{\lambda _{a}}}  - \mu _r^* \\
	 &=\sum\limits_{a \in {\cal A}\backslash {{\cal A}^\textup{tr}}} {c_a^r}  + \sum\limits_{a \in {{\cal O}^\textup{tr}}} {\left( {c_a^r - d_p^a\lambda _a^*} \right)}  + \sum\limits_{a \in {{\cal C}^\textup{tr}}} {\left( {c_a^r - d_p^a{\lambda _a}} \right)}  - \mu _r^*
	\end{split}
	\label{rc:strengthen}
\end{equation}
\endgroup

If only the additional cut is added to BMP, the initial DBSP can be further simplified to obtain the optimal $\lambda_a^*, a\in {\cal O}^\textup{tr}$ and $\mu_r^*, r\in{\cal R}$ faster. Let $\lambda_a, a\in {\cal C}^\textup{tr}$ take $-\infty$, then the constraints in DBSP associated with the paths containing at least one closed arc always remain feasible, and therefore such constraints can be removed. This technique does not influence the optimal solution of DBSP since the coefficients of variables associated with the closed arcs are zeros. When solving PBSP instead of DBSP, such technique implies that we can omit path variables that contain at least one closed arc. This can be realized by forbidding the closed arcs in the passenger routing networks when using column generation. After obtaining the optimal dual variables, $\lambda_a, a\in {\cal C}^\textup{tr}$ can be obtained by solving the cut strengthening model.

\section{Case description}
\subsection{Case Toy}
	\label{app: case toy}
	The topology of the Toy instance is shown in Figure \ref{fig: toy net}. Stations 1, 2, 4, 5, 7\oxc and 8 are set to be terminal stations, and stations 3 and 6 are intermediate ones.
	\begin{figure}[htb]
		\centering
		\FIGURE
		{\includegraphics[width=0.3\linewidth]{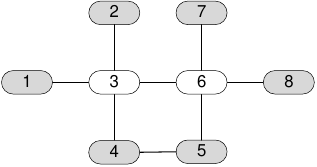}}
		{\centering{Topology of the Toy network}\label{fig: toy net}}
		{}
	\end{figure}

	The cycle time is set as 20 minutes and a total of four periods is considered. After the line planning and timetabling for the single period using existing models in LinTim, eight original trains are operated in each period. We consider the time discretization of of two minutes. The minimum train headway time for the skip-departure case $h_\textup{pd}$ is set as 2 minutes, and the rest are all set as 4 minutes. We set the distance of each section to be 50 kilometers and the running times are set based on the basic setting of the Toy dataset. The additional section running times due to departures and arrivals are both set as 2 minutes, and the transition time of the rolling stock units is set as 6 minutes. 168 time-dependent passenger groups are considered and each group has a preferred departure time period and can accept services within the preferred period, or the half period just before or next to the preferred one. The minimum transfer walking time is set as 4 minutes. As the second period has the highest passenger demand, we set the preferred departure times of extra trains at the center of the second period. Besides, we include all the possible station paths for extra trains.

\subsection{Case Lower Saxony}
\label{app: lower saxony}

The Lower Saxony case in LinTim is designed based on a reginal railway network in Germany, whose topology is shown in Figure \ref{fig: ls net}. Stations 1, 4, 15, 19, 21, 22\oxc and 30 are set as terminal stations where trains can originate or terminate.
\begin{figure}[htb]
	\centering
	\FIGURE
	{\includegraphics[width=\linewidth]{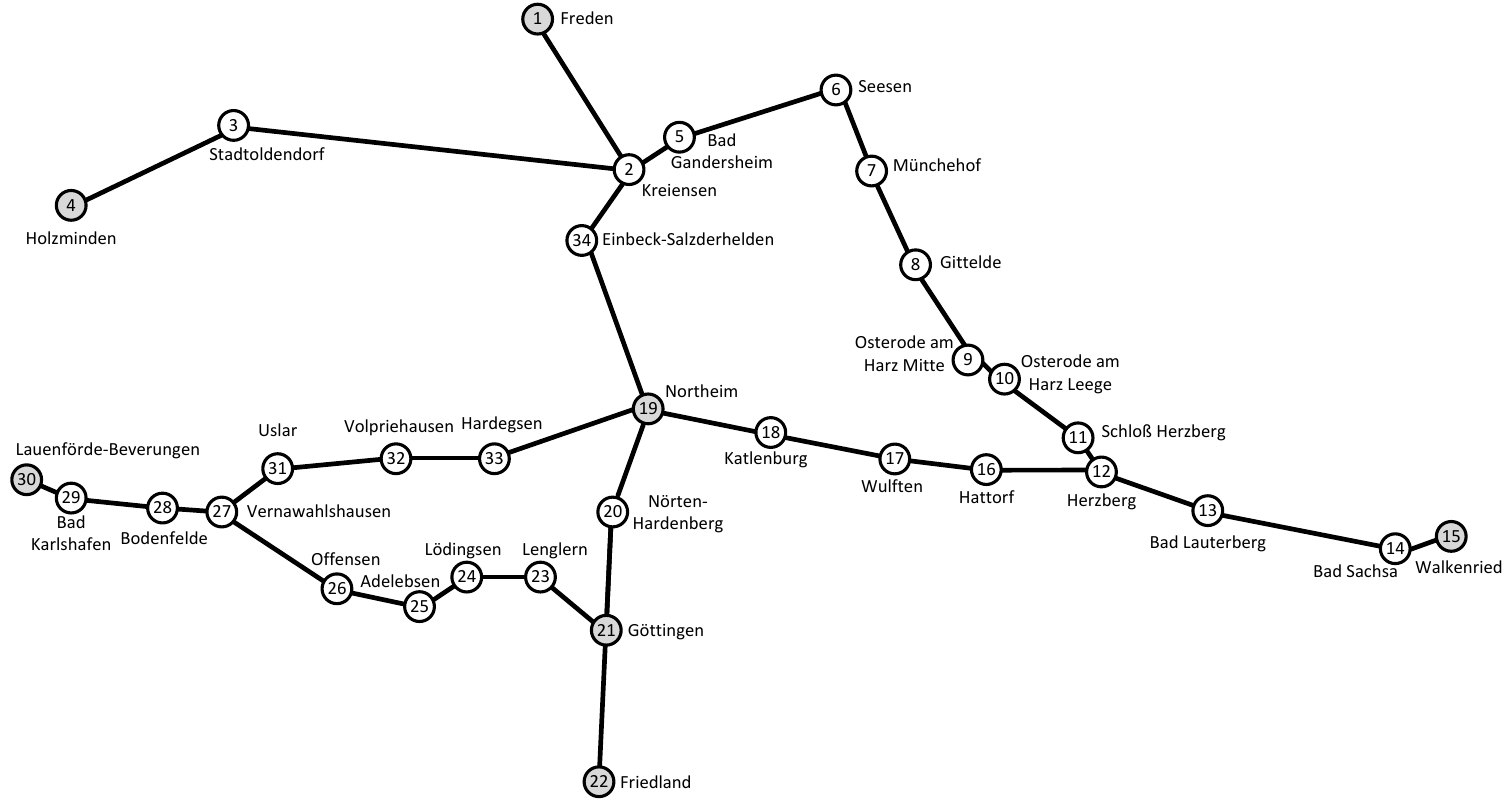}}
	{\centering{Topology of the Lower Saxony network}\label{fig: ls net}}
	{}
\end{figure}

The transition time for rolling stock units is set as 10 minutes. Passengers' preferred and allowable departure time windows are set at one hour and two hours in length, respectively. To generate the time-dependent irregular passenger demand, we first define the peak and non-peak periods, then scale the basic hourly demand from LinTim using different scaling parameters for each period. The resulting time-dependent demand is then adjusted slightly, either increased or decreased, in a random manner. Since the second and third periods are assumed to have the highest passenger demands, the preferred departure times for extra train services are set at the centers of these periods. Other parameter settings remain the same as in the Toy case.

\section{Arc-based formulation}
\label{app: arc formulation}
		In this appendix we give the arc-based formulation for HPTTP. We introduce the passenger arc variable $z_a^r$ which denotes the number of passengers in group $r$ that travel on arc $a$. The arc-based formulation states as follows:
		\begin{align}
			\min \quad&\sum\limits_{r \in {\cal R}}{\sum _{a \in {{\cal A}_r}}} c _a^rz_a^r + \sum\limits_{r \in {\cal R}} {{f_r}{q_r}} \label{obj-arc}\\
			\textrm{s.t.} \quad & \text{(\ref{con:theta x})--(\ref{con:sta inv})}\nonumber\\
			&\sum\limits_{a \in \delta _r^ + (n,t)} {z_a^r}  - \sum\limits_{a \in \delta _r^ - (n,t)} {z_a^r}  = 0 \quad \forall (n,t) \in {\cal {V}_r} \setminus \{ ({o_r},0),({d_r},0)\}
			\label{con: flow bal-arc}\\
			&\sum\limits_{a \in \delta _r^ + ({o_r},0)} {z_a^r + {q_r}}  = {g_r}\quad \forall r \in {\cal R}
			\label{con:pax demand-arc}\\
			&\sum\limits_{{r \in {\cal R}}:a \in {{\cal A}_r}} {z_a^r}  - \sum\limits_{k \in {\cal K}:a \in {{\cal A}_k}} {\sum\limits_{u \in {\cal U}} {{p_u}x_a^{k,u}} }  - \sum\limits_{u \in {\cal U}} {{p_u}y_a^u}  \le 0\quad \forall a \in {{\cal A}^\textup{ex}}
			\label{con:seat cap-xy-arc}\\
			&\sum\limits_{{r \in {\cal R}}:a \in {{\cal A}_r}} {z_a^r}  - \sum\limits_{k \in {\cal K}:a \in {{\cal A}_k}} {\sum\limits_{u \in {\cal U}} {{p_u}x_a^{k,u}} }  \le 0\quad \forall a \in {{\cal A}^\textup{tr}} \backslash {{\cal A}^\textup{ex}}
			\label{con:seat cap-only x-arc}\\
			&\text{(\ref{var x})--(\ref{var w})}, (\ref{var q}) \nonumber\\
			&z^r_a \geq 0 \quad \forall r\in{\cal R}, a \in {\cal A}_r\label{var z}
		\end{align}
		The objective function (\ref{obj-arc}) now uses arc costs to account for passenger travel costs. The passenger flow balance constraints at vertices (\ref{con: flow bal-arc}) are now introduced to represent the path structures of passenger routes. New passenger demand constraints (\ref{con:pax demand-arc}) and seat capacity constraints (\ref{con:seat cap-xy-arc})--(\ref{con:seat cap-only x-arc}) are now using the passenger arc variables. 
\end{APPENDICES}

\end{document}